\documentclass[final,runningheads]{svjour2}
\smartqed  
\usepackage{graphicx}
\usepackage{biograph}        
\journalname{Math. Program., Ser. A}
\usepackage{mathptmx}
\usepackage{subfigure}

\usepackage[usenames]{color}

\usepackage[left=1in,top=1in,right=1in,bottom=1in,letterpaper]{geometry}
\usepackage{setspace,url}
\usepackage[authoryear,square,sort,numbers]{natbib}
\usepackage{amsmath}
\usepackage{amsxtra, amsfonts,amscd, amssymb, graphicx,epsf}
\usepackage{algorithm,algorithmic}

\numberwithin{equation}{section}

\newcommand{\be}{\begin{equation}}
\newcommand{\ee}{\end{equation}}
\newcommand{\ba}{\begin{array}}
\newcommand{\ea}{\end{array}}
\newcommand{\bea}{\begin{eqnarray}}
\newcommand{\eea}{\end{eqnarray}}
\newcommand{\beaa}{\begin{eqnarray*}}
\newcommand{\eeaa}{\end{eqnarray*}}

\newcommand{\half}{\frac{1}{2}}
\newcommand{\br}{\mathbb{R}}

\usepackage[usenames]{color}

\newcounter{alg-counter} 
\newtheorem{algorithm2}[alg-counter]{Algorithm}

 \newcounter{alg-item}

 \newcounter{sub-alg-item}

\usepackage{multirow}
\newcommand{\bean}{\begin{eqnarray}\nonumber}

\newcommand{\dis}{\displaystyle}
\newcommand{\df}[2]{\dis{\frac{#1}{#2}}}
\newcommand{\A}{\mathcal{A}}
\newcommand{\X}{\mathcal{X}}
\newcommand{\zero}{\mathbf{0}}
\newcommand{\Diag}{\mbox{Diag}}
\newcommand{\vvec}{\mbox{vec}}
\newcommand{\Tr}{\mbox{Tr}}
\newcommand{\etal}{{ et al. }}
\newcommand{\argmin}{\mbox{argmin}}

\newcommand{\st}{\mbox{s.t.}}
\newcommand{\mtn}{{m\times n}}
\newcommand{\sgn}{\mbox{sgn}}
\newcommand{\SGN}{\mbox{SGN}}
\newcommand{\rank}{\mbox{rank}}
\begin{document}
\title{Fixed point and Bregman iterative methods for matrix rank minimization}
\titlerunning{Fixed point and Bregman iterative methods for matrix rank minimization}
\author{Shiqian Ma  \and Donald Goldfarb \and Lifeng Chen}
\institute{Department of Industrial Engineering and Operations
Research, Columbia University, New York, NY 10027. Email: \{sm2756,
goldfarb, lc2161\}@columbia.edu\\ Research supported in part by NSF
Grant DMS 06-06712, ONR Grants N00014-03-0514 and N00014-08-1-1118,
and DOE Grants DE-FG01-92ER-25126 and DE-FG02-08ER-58562.}
\date{\bf Submitted October 27, 2008,  Revised May 7, 2009}
\maketitle

\begin{abstract}The linearly constrained matrix rank minimization problem is widely applicable in many fields such as control,
signal processing and system identification. The tightest convex
relaxation of this problem is the linearly constrained nuclear norm
minimization. Although the latter can be cast as a semidefinite
programming problem, such an approach is computationally expensive
to solve when the matrices are large.
In this paper, we propose fixed point and Bregman iterative
algorithms for solving the nuclear norm minimization problem and
prove convergence of the first of these algorithms. By using a
homotopy approach together with an approximate singular value
decomposition procedure, we get a very fast, robust and powerful
algorithm, which we call FPCA (Fixed Point Continuation with
Approximate SVD), that can solve very large matrix rank minimization
problems \footnote{The code can be downloaded from {\em
http://www.columbia.edu/\~{ }sm2756/FPCA.htm} for non-commercial
use.}.
Our numerical results on randomly generated and real matrix
completion problems demonstrate that this algorithm is much faster
and provides much better recoverability than semidefinite
programming solvers such as SDPT3. For example, our algorithm can
recover $1000\times 1000$ matrices of rank 50 with a relative error
of $10^{-5}$ in about 3 minutes by sampling only 20 percent of the
elements. We know of no other method that achieves as good
recoverability. Numerical experiments on online recommendation, DNA
microarray data set and image inpainting problems demonstrate the
effectiveness of our algorithms.

\keywords{Matrix Rank Minimization \and Matrix Completion Problem
\and Nuclear Norm Minimization \and Fixed Point Iterative Method
\and Bregman Distances \and Singular Value Decomposition}
\end{abstract}

\vspace{1cm}\noindent{\bf AMS subject classification.} 65K05, 90C25,
90C06, 93C41, 68Q32

\section{Introduction} The matrix rank minimization problem can be written as
\begin{align*}\ba{ll}\min & \rank(X) \\ s.t. & X \in \mathcal{C},\ea\end{align*}
where $X\in\br^\mtn$ and $\mathcal{C}$ is a convex set. This model
has many applications such as determining a low-order controller for
a plant \citep{ElGhaoui-Gahinet-1993} and a minimum order linear
system realization \citep{Fazel-Hindi-Boyd-2001}, and solving
low-dimensional Euclidean embedding problems
\citep{Linial-London-Rabinovich-1995}.

In this paper, we are interested in methods for solving the affinely
constrained matrix rank minimization problem
\bea\label{prob:min-rank}\ba{cl}
\min & \rank(X)\\
\st & \A(X)=b, \ea\eea where $X\in\br^{m\times n}$ is the decision
variable, and the linear map $\A:\br^{m\times n}\rightarrow\br^p$
and vector $b\in\br^p$ are given.

The matrix completion problem
\bea\label{prob:rank-matrix-completion-intro}\ba{cl}\min & \rank(X)
\\ \st & X_{ij} = M_{ij}, (i,j)\in\Omega \ea\eea is a special case of
\eqref{prob:min-rank}, where $X$ and $M$ are both $m\times n$
matrices and $\Omega$ is a subset of index pairs $(i,j).$ The so
called collaborative filtering problem
\citep{Rennie-Srebro-2005,Srebro-thesis-2004} can be cast as a
matrix completion problem. Suppose users in an online survey provide
ratings of some movies. This yields a matrix $M$ with users as rows
and movies as columns whose $(i,j)$-th entry $M_{ij}$ is the rating
given by the $i$-th user to the $j$-th movie. Since most users rate
only a small portion of the movies, we typically only know a small
subset $\{M_{ij}|(i,j)\in\Omega\}$ of the entries. Based on the
known ratings of a user, we want to predict the user's ratings of
the movies that the user did not rate; i.e., we want to fill in the
missing entries of the matrix. It is commonly believed that only a
few factors contribute to an individual's tastes or preferences for
movies. Thus the rating matrix $M$ is likely to be of {\em
numerical} low rank in the sense that relatively few of the top
singular values account for most of the sum of all of the singular
values. Finding such a low-rank matrix $M$ corresponds to solving
the matrix completion problem
\eqref{prob:rank-matrix-completion-intro}.

{\bf \subsection{Connections to compressed sensing}}

When the matrix $X$ is diagonal, problem \eqref{prob:min-rank}
reduces to the cardinality minimization problem
\bea\label{prob:cardinality-minimization}\ba{cl}\min & \|x\|_0
\\ \st & Ax=b,\ea\eea where $x\in\br^n, A\in\br^{m\times n},
b\in\br^m$ and $\|x\|_0$ denotes the number of nonzeros in the
vector $x.$ This problem finds the sparsest solution to an
underdetermined system of equations and has a wide range of
applications in signal processing. This problem is NP-hard
\citep{Natarajan-1995}. To get a more computationally tractable
problem, we can replace $\|x\|_0$ by its convex envelope.

\begin{definition}\label{definition:convex-envelope}
The convex envelope of a function $f:\mathcal{C}\rightarrow\br$ is
defined as the largest convex function $g$ such that $g(x)\leq f(x)$
for all $x\in\mathcal{C}$ (see e.g.,
\citep{Hiriart-Urruty-Lemarechal-1993}).
\end{definition}

It is well known that the convex envelope of $\|x\|_0$  is
$\|x\|_1$, the $\ell_1$ norm of $x$, which is the sum of the
absolute values of all components of $x$. Replacing the objective
function $\|x\|_0$ in \eqref{prob:cardinality-minimization} by
$\|x\|_1$ yields the so-called basis pursuit problem
\bea\label{prob:basis-pursuit-intro}\ba{cl}\min & \|x\|_1 \\
\st & Ax=b.\ea\eea The basis pursuit problem has received an
increasing amount of attention since the emergence of the field of
compressed sensing (CS) \citep{Candes-Romberg-Tao-2006,Donoho-2006}.
Compressed sensing theories connect the NP-hard problem
\eqref{prob:cardinality-minimization} to the convex and
computationally tractable problem \eqref{prob:basis-pursuit-intro}
and provide guarantees for when an optimal solution to
\eqref{prob:basis-pursuit-intro} gives an optimal solution to
\eqref{prob:cardinality-minimization}. In the cardinality
minimization and basis pursuit problems
\eqref{prob:cardinality-minimization} and
\eqref{prob:basis-pursuit-intro}, $b$ is a vector of measurements of
the signal $x$ obtained using the sampling matrix $A$. The main
result of compressed sensing is that when the signal $x$ is sparse,
i.e., $k:=\|x\|_0\ll n,$ we can recover the signal by solving
\eqref{prob:basis-pursuit-intro} with a very limited number of
measurements, i.e., $m\ll n$, when $A$ is a Gaussian random matrix
or when it corresponds to a partial Fourier transformation. Note
that if $b$ is contaminated by noise, the constraint $Ax=b$ in
\eqref{prob:basis-pursuit-intro} must be relaxed, resulting in
either the problem \bea\label{prob:basis-pursuit-noise}\ba{cl}\min &
\|x\|_1
\\ \st & \|Ax-b\|_2\leq\theta\ea\eea or its Lagrangian
version \bea\label{prob:basis-pursuit-lagrange}\min
\mu\|x\|_1+\half\|Ax-b\|_2^2,\eea where $\theta$ and $\mu$ are
parameters and $\|x\|_2$ denotes the Euclidean norm of a vector
$x$.. Algorithms for solving \eqref{prob:basis-pursuit-intro} and
its variants \eqref{prob:basis-pursuit-noise} and
\eqref{prob:basis-pursuit-lagrange} have been widely investigated
and many algorithms have been suggested including convex
optimization methods
(\citep{Candes-Romberg-2005-l1-magic,Figueiredo-Nowak-Wright-2007,Hale-Yin-Zhang-2007,Kim-Koh-Lustig-Boyd-Gorinevsky-2007,vandenBerg-Friedlander-2008})
and heuristic methods
(\citep{Tibshirani-1996,Donoho-Tsaig-Drori-Starck-2006,Tropp-2006,Donoho-Tsaig-2006,Dai-Milenkovie-2008}).

{\bf \subsection{Nuclear norm minimization}}

The rank of a matrix is the number of its positive singular values.
The matrix rank minimization \eqref{prob:min-rank} is NP-hard in
general due to the combinational nature of the function
$\rank(\cdot)$. Similar to the cardinality function $\|x\|_0$, we
can replace $\rank(X)$ by its convex envelope to get a convex and
more computationally tractable approximation to
\eqref{prob:min-rank}. It turns out that the convex envelope of
$\rank(X)$ on the set $\{X\in\br^\mtn: \|X\|_2\leq 1\}$ is the
nuclear norm $\|X\|_*$ \citep{Fazel-thesis-2002}, i.e., the nuclear
norm is the best convex approximation of the rank function over the
unit ball of matrices with norm less than one, where $\|X\|_2$ is
the operator norm of $X$. The nuclear norm and operator norm are
defined as follows.

\begin{definition}\label{definition:nuclear-norm-operator-norm}
Nuclear norm and Operator norm. Assume that the matrix $X$ has $r$
positive singular values of
$\sigma_1\geq\sigma_2\geq\ldots\geq\sigma_r>0$. The nuclear norm of
$X$ is defined as the sum of its singular values, i.e., \bean
\|X\|_*:=\sum_{i=1}^r\sigma_i(X).\eea The operator norm of matrix
$X$ is defined as the largest singular value of $X$, i.e.,
\beaa\|X\|_2:=\sigma_1(X).\eeaa
\end{definition}

The nuclear norm is also known as Schatten 1-norm or Ky Fan norm.
Using it as an approximation to $\rank(X)$ in \eqref{prob:min-rank}
yields the nuclear norm minimization problem
\bea\label{prob:nuclear-norm-minimization-intro}\ba{cl}\min &
\|X\|_* \\ \st & \A(X)=b.\ea\eea As in the basis pursuit problem, if
$b$ is contaminated by noise, the constraint $\A(X)=b$ must be
relaxed, resulting in either the problem \bean\ba{cl}\min & \|X\|_*
\\ \st & \|\A(X)-b\|_2\leq\theta\ea\eea or its Lagrangian version
\bea\label{prob:min-unconstrained}\min
\mu\|X\|_*+\df{1}{2}\|\A(X)-b\|_2^2,\eea where $\theta$ and $\mu$
are parameters.

Note that if we write $X$ in vector form by stacking the columns of
$X$ in a single vector $\vvec(X)\in\br^{mn}$, then we get the
following equivalent formation of
\eqref{prob:nuclear-norm-minimization-intro}:
\bea\label{prob:nuclear-norm-vecX}\ba{cl}\min & \|X\|_*
\\\st & A\ \vvec(X)=b, \ea\eea where $A\in\br^{p\times
mn}$ is the matrix corresponding to the linear map $\A.$ An
important question is: when will an optimal solution to the nuclear
norm minimization problem
\eqref{prob:nuclear-norm-minimization-intro} give an optimal
solution to matrix rank minimization problem \eqref{prob:min-rank}.
In response to this question, Recht \etal
\citep{Recht-Fazel-Parrilo-2007} proved that if the entries of $A$
are suitably random, e.g., i.i.d. Gaussian, then with very high
probability, most $m\times n$ matrices of rank $r$ can be recovered
by solving the nuclear norm minimization
\eqref{prob:nuclear-norm-minimization-intro} or equivalently,
\eqref{prob:nuclear-norm-vecX}, whenever $p\geq Cr(m+n)\log(mn),$
where $C$ is a positive constant.

For the matrix completion problem
\eqref{prob:rank-matrix-completion-intro}, the corresponding nuclear
norm minimization problem is
\bea\label{prob:nuclear-norm-matrix-completion-intro}\ba{cl}\min &
\|X\|_* \\ \st & X_{ij} = M_{ij}, (i,j) \in \Omega. \ea\eea Cand\`es
\etal \citep{Candes-Recht-2008} proved the following result.
\begin{theorem}\label{thm:Candes-Recht}
Let $M$ be an $n_1\times n_2$ matrix of rank $r$ with SVD
$$M=\sum_{k=1}^r\sigma_ku_kv_k^\top,$$ where the family $\{u_k\}_{1\leq k\leq r}$
is selected uniformly at random among all families of $r$
orthonormal vectors, and similarly for the family $\{v_k\}_{1\leq
k\leq r}$. Let $n=\max(n_1,n_2)$. Suppose we observe $m$ entries of
$M$ with locations sampled uniformly at random. Then there are
constants $C$ and $c$ such that if
$$m\geq Cn^{5/4}r\log n,$$ the minimizer to the problem
\eqref{prob:nuclear-norm-matrix-completion-intro} is unique and
equal to $M$ with probability at least $1-cn^{-3}$. In addition, if
$r\leq n^{1/5}$, then the recovery is exact with probability at
least $1-cn^{-3}$ provided that $$m\geq Cn^{6/5}r\log n.$$
\end{theorem}
This theorem states that a surprisingly small number of entries are
sufficient to complete a low-rank matrix with high probability.

Recently, this result was strengthened by Cand\`es and Tao in
\citep{Candes-Tao-2009}, where it is proved that under certain
incoherence conditions, the number of samples $m$ that are required
is only $O(nr\log n).$


The dual problem corresponding to the nuclear norm minimization
problem \eqref{prob:nuclear-norm-minimization-intro} is
\bea\label{prob:dual-nuclear-norm-minimization}\ba{cl} \max & b^\top
z \\ \st & \|\A^*(z)\|_2\leq 1, \ea\eea where $\A^*$ is the adjoint
operator of $\A$. Both \eqref{prob:nuclear-norm-minimization-intro}
and \eqref{prob:dual-nuclear-norm-minimization} can be rewritten as
equivalent semidefinite programming (SDP) problems. The SDP
formulation of \eqref{prob:nuclear-norm-minimization-intro} is:
\bea\label{prob:sdp-nuclear-norm-primal}\ba{cl}
\displaystyle\min_{X,W_1,W_2}& \half(\Tr(W_1)+\Tr(W_2)) \\\st &
\begin{bmatrix}W_1 & X \\ X^\top & W_2\end{bmatrix}\succeq 0 \\ & \A(X)=b,
\ea\eea where $\Tr(X)$ denotes the trace of the square matrix $X$.
The SDP formulation of \eqref{prob:dual-nuclear-norm-minimization}
is: \bea\label{prob:sdp-nuclear-norm-dual}\ba{cl}\displaystyle\max_z
& b^\top z
\\ \st & \begin{bmatrix} I_m & \A^*(z) \\ \A^*(z)^\top & I_n \end{bmatrix} \succeq 0. \ea\eea

Thus to solve \eqref{prob:sdp-nuclear-norm-primal} and
\eqref{prob:sdp-nuclear-norm-dual}, we can use SDP solvers such as
SeDuMi \citep{Sturm-1999} and SDPT3 \citep{Tutuncu-Toh-Todd-2003} to
solve \eqref{prob:sdp-nuclear-norm-primal} and
\eqref{prob:sdp-nuclear-norm-dual}. Note that the number of
variables in \eqref{prob:sdp-nuclear-norm-primal} is
$\half(m+n)(m+n+1)$. SDP solvers cannot usually solve a problem when
$m$ and $n$ are both much larger than $100.$

Recently, Liu and Vandenberghe \citep{Liu-Vandenberghe-2008}
proposed an interior-point method for another nuclear norm
approximation problem \bea\label{prob:vandenberghe-nuclear-norm}\min
\|\A(x)-B\|_*,\eea where $B\in\br^\mtn$ and
$$\A(x)=x_1A_1+x_2A_2+\cdots+x_pA_p$$ is a linear mapping from
$\br^p$ to $\br^\mtn.$ The equivalent SDP formulation of
\eqref{prob:vandenberghe-nuclear-norm} is
\bea\label{prob:sdp-vandenberghe-nuclear-norm}\ba{cl}\displaystyle\min_{x,W_1,W_2}
& \half(\Tr(W_1)+\Tr(W_2)) \\ \st & \begin{bmatrix} W_1 &
(\A(x)-B)^\top \\ \A(x)-B & W_2
\end{bmatrix}\succeq 0. \ea\eea
Liu and Vandenberghe \citep{Liu-Vandenberghe-2008} proposed a
customized method for computing the scaling direction in an interior
point method for solving the SDP
\eqref{prob:sdp-vandenberghe-nuclear-norm}. The complexity of each
iteration in their method was reduced from $O(p^6)$ to $O(p^4)$ when
$m=O(p)$ and $n=O(p)$; thus they were able to solve problems up to
dimension $m=n=350.$ 

Another algorithm for solving
\eqref{prob:nuclear-norm-minimization-intro} is due to Burer and
Monteiro \citep{Burer-Monteiro-2003, Burer-Monteiro-2005}, (see also
Rennie and Srebro \citep{Rennie-Srebro-2005, Srebro-thesis-2004}).
This algorithm uses the low-rank factorization $X=LR^\top$ of the
matrix $X\in\br^\mtn$, where $L\in\br^{m\times r}, R\in\br^{n\times
r},r\leq\min\{m,n\},$ and solves the optimization problem
\bea\label{prob:low-rank-factorization}\ba{cl}\displaystyle\min_{L,R}&
\half(\|L\|_F^2+\|R\|_F^2) \\ \st & \A(LR^\top)=b, \ea\eea where
$\|X\|_F$ denotes the Frobenius norm of the matrix $X$: \beaa\|X\|_F
:=
(\sum_{i=1}^r\sigma_i^2)^{1/2}=(\sum_{i,j}X_{ij}^2)^{1/2}=(\Tr(XX^\top))^{1/2}.\eeaa

It is known that as long as $r$ is chosen to be sufficiently larger
than the rank of the optimal solution matrix of the nuclear norm
problem \eqref{prob:nuclear-norm-minimization-intro}, this low-rank
factorization problem is equivalent to the nuclear norm problem
\eqref{prob:nuclear-norm-minimization-intro} (see e.g.,
\citep{Recht-Fazel-Parrilo-2007}). The advantage of this low-rank
factorization formulation is that both the objective function and
the constraints are differentiable. Thus gradient-based optimization
algorithms such as conjugate gradient algorithms and augmented
Lagrangian algorithms can be used to solve this problem. However,
the constraints in this problem are nonconvex, so one can only be
assured of obtaining a local minimizer. Also, how to choose $r$ is
still an open question.

One very interesting algorithm is the so called singular value
thresholding algorithm (SVT) \citep{Cai-Candes-Shen-2008} which
appeared almost simultaneously with our work. SVT is inspired by the
linearized Bregman algorithms for compressed sensing and
$\ell_1$-regularized problems. In \citep{Cai-Candes-Shen-2008} it is
shown that SVT is efficient for large matrix completion problems.
However, SVT only works well for very low rank matrix completion
problems. For problems where the matrices are not of very low rank,
SVT is slow and not robust therefore often fails.

Our algorithms have some similarity with the SVT algorithm in that
they make use of {\em matrix shrinkage} (see Section 2). However,
other than that, they are greatly different. All of our methods are
based on a fixed point continuation (FPC) algorithm which uses an
operator splitting technique for solving
\eqref{prob:min-unconstrained}. By adopting a Monte Carlo
approximate SVD in the FPC, we get an algorithm, which we call FPCA
(Fixed Point Continuation with Approximate SVD), that usually gets
the optimal solution to \eqref{prob:min-rank} even if the condition
of Theorem \ref{thm:Candes-Recht}, or those for the affine
constrained case, are violated. Moreover, our algorithm is much
faster than state-of-the-art SDP solvers such as SDPT3 applied to
\eqref{prob:sdp-nuclear-norm-primal}. Also, FPCA can recover
matrices of moderate rank that cannot be recovered by SDPT3, SVT,
etc. with the same amount of samples. For example, for matrices of
size $1000\times 1000$ and rank 50, FPCA can recover them with a
relative error of $10^{-5}$ in about 3 minutes by sampling only 20
percent of the matrix elements. As far as we know, there is no other
method that has as good a recoverability property.

{\bf \subsection{Outline and Notation}}

{\bf Outline.} The rest of this paper is organized as follows. In
Section 2 we review the fixed-point continuation algorithm for
$\ell_1$-regularized problems. In Section 3 we give an analogous
fixed-point iterative algorithm for the nuclear norm minimization
problem and prove that it converges to an optimal solution. In
Section 4 we discuss a continuation technique for accelerating the
convergence of our algorithm. In Section 5 we propose a Bregman
iterative algorithm for nuclear norm minimization extending the
approach in \citep{Yin-Osher-Goldfarb-Darbon-2008} for compressed
sensing to the rank minimization problem. In Section 6 we
incorporate a Monte-Carlo approximate SVD procedure into our
fixed-point continuation algorithm to speed it up and improve its
ability to recover low-rank matrices. Numerical results for both
synthesized matrices and real problems are given in Section 7. We give conclusions in Section 8.  

{\bf Notation.} 
Throughout this paper, we always assume that the singular values are
arranged in nonincreasing order, i.e.,
$\sigma_1\geq\sigma_2\geq\ldots\geq\sigma_r>0=\sigma_{r+1}=\ldots=\sigma_{\min\{m,n\}}.$
$\partial f$ denotes the subdifferential of the function $f$ and
$g^k = g(X^k) = \A^*(\A(X^k)-b))$ is the gradient of function
$\half\|\A(X)-b\|_2^2$ at the point $X^k$. $\Diag(s)$ denotes the
diagonal matrix whose diagonal elements are the elements of the
vector $s$. $\sgn(t)$ is the signum function of $t\in\br$, i.e.,
\beaa\sgn(t):=\left\{\ba{ll} +1 & \mbox{ if } t
> 0, \\ 0 & \mbox{ if } t =0, \\ -1 & \mbox{ if } t<0,
\ea\right.\eeaa while the signum multifunction of $t\in\br$ is
\beaa\SGN(t):=\partial |t| = \left\{\ba{ll} \{+1\} & \mbox{ if } t >
0, \\\ [-1,1] & \mbox{ if } t =0, \\ \{-1\} & \mbox{ if } t<0.
\ea\right.\eeaa We use $a\odot b$ to denote the elementwise
multiplication of two vectors $a$ and $b$. We use $X(k:l)$ to denote
the submatrix of $X$ consisting of the $k$-th to $l$-th column of
$X$. We use $\br_+^n$ to denote the nonnegative orthant in $\br^n.$

\section{Fixed point iterative algorithm}
Our fixed point iterative algorithm for solving
\eqref{prob:min-unconstrained} is the following simple two-line
algorithm:
\bea\label{eq:fpc-one-step-scheme}\left\{\ba{l}Y^k=X^k-\tau
g(X^k)\\X^{k+1}=S_{\tau\mu}(Y^k),\ea\right.\eea where
$S_{\nu}(\cdot)$ is the matrix shrinkage operator which will be
defined later.

Our algorithm \eqref{eq:fpc-one-step-scheme} is inspired by the
fixed point iterative algorithm proposed in
\citep{Hale-Yin-Zhang-2007} for the $\ell_1$-regularized problem
\eqref{prob:basis-pursuit-lagrange}. The idea behind this algorithm
is an operator splitting technique. Note that $x^*$ is an optimal
solution to \eqref{prob:basis-pursuit-lagrange} if and only if
\bea\label{optcond:ell_1-regularized}\zero\in\mu\SGN(x^*)+g^*,\eea
where $g^*=A^\top(Ax^*-b)$. For any $\tau>0$,
\eqref{optcond:ell_1-regularized} is equivalent to
\bea\label{optcond:ell_1-regularized-tau}\zero\in\tau\mu\SGN(x^*)+\tau
g(x^*). \eea Note that the operator
$T(\cdot):=\tau\mu\SGN(\cdot)+\tau g(\cdot)$ on the right hand side
of \eqref{optcond:ell_1-regularized-tau} can be split into two
parts: $T(\cdot)=T_1(\cdot)-T_2(\cdot),$ where
$T_1(\cdot)=\tau\mu\SGN(\cdot)+I(\cdot)$ and
$T_2(\cdot)=I(\cdot)-\tau g(\cdot)$.

Letting $y=T_2(x^*)=x^*-\tau A^\top(Ax^*-b)$,
\eqref{optcond:ell_1-regularized-tau} is equivalent to
\bea\label{optcond:ell_1-regularized-y-x-taug}\zero\in
T_1(x^*)-y=\tau\mu\SGN(x^*) + x^*-y.\eea Note that
\eqref{optcond:ell_1-regularized-y-x-taug} is actually the
optimality conditions for the following convex problem
\bea\label{prob:ell_1-regularized-shrinkage}\min_{x^*}
\tau\mu\|x^*\|_1+\half\|x^*-y\|_2^2.\eea This problem has a closed
form optimal solution given by the so called shrinkage operator:
\beaa x^*=\tilde{s}_\nu(y),\eeaa where $\nu=\tau\mu,$ and shrinkage
operator $\tilde{s}_\nu(\cdot)$ is given by
\bea\label{eq:shrinkage-ell-1}\tilde{s}_\nu(\cdot)=\sgn(\cdot)\odot\max\{|\cdot|-\nu,0\}.\eea
Thus, the fixed point iterative algorithm is given by
\bea\label{iterative-fixed-point-ell-1}
x^{k+1}=\tilde{s}_{\tau\mu}(x^k-\tau g^k). \eea Hale \etal
\citep{Hale-Yin-Zhang-2007} proved global and finite convergence of
this algorithm to the optimal solution of the $\ell_1$-regularized
problem \eqref{prob:basis-pursuit-lagrange}.

Motivated by this work, we develop a fixed point iterative algorithm
for \eqref{prob:min-unconstrained}. Since the objective function in
\eqref{prob:min-unconstrained} is convex, $X^*$ is the optimal
solution to \eqref{prob:min-unconstrained} if and only if
\bea\label{OPTcond:convex-nuclear}\zero\in\mu\partial\|X^*\|_*+g(X^*),\eea
where $g(X^*)=\A^*(\A(X^*)-b)$. Note that if the Singular Value
Decomposition (SVD) of $X$ is $X=U\Sigma V^\top$, where
$U\in\br^{m\times r},\Sigma=\Diag(\sigma)\in\br^{r\times r},
V\in\br^{n\times r},$ then (see e.g.,
\citep{Borwein-Lewis-2000-book,Bach-2008})
\bean\partial\|X\|_*=\{UV^\top+W:U^\top W=0,WV=0,\|W\|\leq 1\}. \eea
Hence, we get the following optimality conditions for
\eqref{prob:min-unconstrained}:
\begin{theorem}\label{thm:OPTCondThm2}
The matrix $X\in\br^\mtn$ with singular value decomposition
$X=U\Sigma V^\top,$ $U\in\br^{m\times
r},\Sigma=\Diag(\sigma)\in\br^{r\times r}, V\in\br^{n\times r},$ is
optimal for the problem \eqref{prob:min-unconstrained} if and only
if there exists a matrix $W\in\br^\mtn$ such that
\begin{subequations}
\begin{align}
\mu(UV^\top+W)+g(X)=0,\label{optcond:exist-W-unconstrained-problem-a}
\\  U^\top W=0,WV=0,\|W\|_2\leq
1.\label{optcond:exist-W-unconstrained-problem-b}
\end{align}
\end{subequations}
\end{theorem}

Now based on the optimality conditions
\eqref{OPTcond:convex-nuclear}, we can develop a fixed point
iterative scheme for solving \eqref{prob:min-unconstrained} by
adopting the operator splitting technique described at the beginning
of this section. Note that \eqref{OPTcond:convex-nuclear} is
equivalent to
\bea\label{optcond:split}\zero\in\tau\mu\partial\|X^*\|_*
+X^*-(X^*-\tau g(X^*))\eea for any $\tau > 0$. If we let \bean Y^* =
X^*-\tau g(X^*),\eea then \eqref{optcond:split} is reduced to
\bea\label{optcond:reduced}\zero\in\tau\mu\partial\|X^*\|_*
+X^*-Y^*,\eea i.e., $X^*$ is the optimal solution to
\bea\label{prob:unique-solution-shrinkage-Y*}\min_{X\in\br^{m\times
n}}\tau\mu\|X\|_*+\df{1}{2}\|X-Y^*\|_F^2\eea

In the following we will prove that the matrix shrinkage operator
applied to $Y^*$ gives the optimal solution to
\eqref{prob:unique-solution-shrinkage-Y*}. First, we need the
following definitions.
\begin{definition}
[\bf Nonnegative Vector Shrinkage Operator] Assume $x\in\br_+^{n}$.
For any $\nu>0$, the nonnegative vector shrinkage operator
$s_\nu(\cdot)$ is defined as
$$s_\nu(x):=\bar{x},\mbox{ with }\bar{x}_i=\left\{\ba{ll}x_i-\nu,&\mbox{ if }x_i-\nu>0 \\0,&\mbox{ o.w. }\ea\right.$$
\end{definition}

\begin{definition}
[\bf Matrix Shrinkage Operator] Assume $X\in\br^{m\times n}$ and the
SVD of $X$ is given by $X=U\Diag(\sigma)V^\top$, $U\in\br^{m\times
r},\sigma\in\br_+^r,V\in\br^{n\times r}.$ For any $\nu>0$, the
matrix shrinkage operator $S_\nu(\cdot)$ is defined as \bean
S_\nu(X):=U\Diag(\bar\sigma) V^\top,\qquad\mbox{with
}\bar\sigma=s_\nu(\sigma).\eea
\end{definition}

\begin{theorem}
Given a matrix $Y\in\br^{m\times n}$ with $\rank(Y)=t$, let its
Singular Value Decomposition (SVD) be $Y=U_Y\Diag(\gamma) V_Y^\top$,
where $U_Y\in\br^{m\times t},\gamma\in\br_+^t,V_Y\in\br^{n\times
t}$, and a scalar $\nu>0$. Then \bea\label{opt:X-optimal-shrinkage}X
:= S_\nu(Y) = U_Y\Diag(s_\nu(\gamma))V_Y^\top\eea is an optimal
solution of the problem
\bea\label{prob:shrinkage-solution}\min_{X\in\br^{m\times
n}}f(X):=\nu\|X\|_*+\df{1}{2}\|X-Y\|_F^2.\eea
\end{theorem}
\proof Without loss of generality, we assume $m\leq n$. Suppose that
the solution $X\in\br^{m\times n}$ to problem
\eqref{prob:shrinkage-solution} has the SVD $X=U\Diag(\sigma)
V^\top,$ where $U\in\br^{m\times
r},\sigma\in\br_+^r,V\in\br^{n\times r}$. Hence, $X$ must satisfy
the optimality conditions for \eqref{prob:shrinkage-solution} which
are \bean\zero\in\nu\partial\|X\|_*+X-Y;\eea i.e., there exists a
matrix
\[W=\bar{U}\begin{bmatrix}\Diag(\bar{\sigma}) & 0 \end{bmatrix}\bar{V}^\top,\] where
$\bar{U}\in\br^{m\times(m-r)}, \bar{V}\in\br^{n\times(n-r)},
\bar{\sigma}\in\br_+^{m-r}$, $\|\bar{\sigma}\|_\infty\leq 1$ and
both $\hat{U}=[U,\bar{U}]$ and $\hat{V}=[V,\bar{V}]$ are orthogonal
matrices,
such that \bea\label{optcond:shrinkage-unique-optcond}\zero=\nu(U
V^\top+W)+X-Y.\eea Hence,
\bea\label{eq:proof-shrinkage-solution-svd-Y}
\hat{U}\begin{bmatrix}\nu I+\Diag(\sigma)&0
&0\\0&\nu\Diag(\bar{\sigma})&0\end{bmatrix}\hat{V}^\top -
U_Y\Diag(\gamma)V_Y^\top = \zero.\eea

To verify that \eqref{opt:X-optimal-shrinkage} satisfies
\eqref{eq:proof-shrinkage-solution-svd-Y}, consider the following
two cases:

{\bf Case 1:} $\gamma_1\geq\gamma_2\geq\ldots\geq\gamma_t>\nu.$ In
this case, choosing $X$ as above, with $r=t,U=U_Y,V=V_Y$ and
$\sigma=s_\nu(\gamma)=\gamma-\nu e$, where $e$ is a vector of $r$
ones, and choosing $\bar\sigma=0$ (i.e., $W=\zero$) satisfies
\eqref{eq:proof-shrinkage-solution-svd-Y}.

{\bf Case 2:}
$\gamma_1\geq\gamma_2\geq\ldots\geq\gamma_k>\nu\geq\gamma_{k+1}\geq\ldots\geq\gamma_t.$
In this case, by choosing $r=k,\hat U(1:t)=U_Y,\hat
V(1:t)=V_Y,\sigma=s_\nu((\gamma_1,\ldots,\gamma_k))$ and
$\bar\sigma_1=\gamma_{k+1}/\nu,\ldots,\bar\sigma_{t-k}=\gamma_t/\nu,\bar\sigma_{t-k+1}=\ldots=\bar\sigma_{m-r}=0,$
$X$ and $W$ satisfy \eqref{eq:proof-shrinkage-solution-svd-Y}.

Note that in both cases, $X$ can be written as the form in
\eqref{opt:X-optimal-shrinkage} based on the way we construct $X$.
\qed 

Based on the above we obtain the fixed point iterative scheme
\eqref{eq:fpc-one-step-scheme} stated at the beginning of this
section for solving problem \eqref{prob:min-unconstrained}.

Moreover, from the discussion following Theorem
\ref{thm:OPTCondThm2} we have
\begin{corollary}\label{corollary:1}
$X^*$ is an optimal solution to problem
\eqref{prob:min-unconstrained} if and only if
$X^*=S_{\tau\mu}(h(X^*))$, where $h(\cdot)=I(\cdot)-\tau g(\cdot).$
\end{corollary}

\section{Convergence results}
In this section, we analyze the convergence properties of the fixed
point iterative scheme \eqref{eq:fpc-one-step-scheme}. Before we
prove the main convergence result, we need some lemmas.
\begin{lemma}\label{lemma:shrinkage non-expansive}
The shrinkage operator $S_\nu$ is non-expansive, i.e., for any $Y_1$
and $Y_2\in\br^\mtn$,
\bea\label{ieq:decrese-shrinkage}\|S_\nu(Y_1)-S_\nu(Y_2)\|_F\leq
\|Y_1-Y_2\|_F.\eea Moreover, 
\bea\label{eq:shrinkage-no-shrinkage}\|Y_1-Y_2\|_F=\|S_\nu(Y_1)-S_\nu(Y_2)\|_F\Longleftrightarrow
Y_1-Y_2=S_\nu(Y_1)-S_\nu(Y_2).\eea
\end{lemma}
\proof Without loss of generality, we assume $m\leq n$. Assume SVDs
of $Y_1$ and $Y_2$ are $Y_1=U_1\Sigma V_1^\top$ and $Y_2=U_2\Gamma
V_2^\top$, respectively, where
$$\Sigma
=\begin{pmatrix} \Diag(\sigma) & 0 \\ 0 & 0
\end{pmatrix}\in\br^\mtn,\Gamma=\begin{pmatrix}\Diag(\gamma)&0\\0&0\end{pmatrix}\in\br^\mtn,$$
$\sigma=(\sigma_1,\ldots,\sigma_s),\sigma_1\geq\ldots\geq\sigma_s>0$
and
$\gamma=(\gamma_1,\ldots,\gamma_t),\gamma_1\geq\ldots\geq\gamma_t>0$.
Note that here $U_1,V_1,U_2$ and $V_2$ are (full) orthogonal
matrices; $\Sigma, \Gamma\in\br^{m\times n}$. Suppose that
$\sigma_1\geq\ldots\geq\sigma_k\geq\nu>\sigma_{k+1}\geq\ldots\geq\sigma_s$
and
$\gamma_1\geq\ldots\geq\gamma_l\geq\nu>\gamma_{l+1}\geq\ldots\geq\gamma_t$,
then
$$\bar{Y}_1:=S_\nu(Y_1)=U_1\bar{\Sigma}V_1^\top,
\bar{Y}_2:=S_\nu(Y_2)=U_2\bar{\Gamma}V_2^\top,$$ where
$$\bar\Sigma
=\begin{pmatrix} \Diag(\bar\sigma) & 0 \\ 0 & 0
\end{pmatrix}\in\br^\mtn,\bar\Gamma=\begin{pmatrix}\Diag(\bar\gamma)&0\\0&0\end{pmatrix}\in\br^\mtn,$$
$\bar\sigma=(\sigma_1-\nu,\ldots,\sigma_k-\nu)$ and $
\bar\gamma=(\gamma_1-\nu,\ldots,\gamma_l-\nu)$. Thus, \beaa\ba{lll}
\|Y_1-Y_2\|_F^2-\|\bar{Y}_1-\bar{Y}_2\|_F^2&=&\Tr((Y_1-Y_2)^\top(Y_1-Y_2))-\Tr((\bar{Y}_1-\bar{Y}_2)^\top(\bar{Y}_1-\bar{Y}_2))\\
&=&\Tr(Y_1^\top Y_1-\bar{Y}_1^\top\bar{Y}_1+Y_2^\top Y_2-\bar{Y}_2^\top\bar{Y}_2)-2\Tr(Y_1^\top Y_2-\bar{Y}_1^\top\bar{Y}_2)\\
&=&\dis\sum_{i=1}^s\sigma_i^2-\sum_{i=1}^k(\sigma_i-\nu)^2+\sum_{i=1}^t\gamma_i^2-\sum_{i=1}^l(\gamma_i-\nu)^2-2\Tr(Y_1^\top
Y_2-\bar{Y}_1^\top\bar{Y}_2)\ea\eeaa

We note that
\begin{equation*}\ba{lll} \Tr(Y_1^\top
Y_2-\bar{Y}_1^\top\bar{Y}_2)&=&\Tr((Y_1-\bar{Y_1})^\top(Y_2-\bar{Y_2})+(Y_1-\bar{Y_1})^\top\bar{Y_2}+\bar{Y_1}^\top(Y_2-\bar{Y_2}))\\
                            &=&\Tr(V_1(\Sigma-\bar\Sigma)^\top U_1^\top U_2(\Gamma-\bar\Gamma)V_2^\top+V_1(\Sigma-\bar\Sigma)^\top U_1^\top U_2\bar\Gamma
                            V_2^\top+V_1{\bar\Sigma}^\top U_1^\top U_2(\Gamma-\bar\Gamma)V_2^\top\\
                            &=&\Tr((\Sigma - \bar\Sigma)^\top U
                            (\Gamma-\bar\Gamma)V^\top + (\Sigma-\bar\Sigma)^\top U\bar\Gamma V^\top + {\bar\Sigma}^\top U(\Gamma-\bar\Gamma)V^\top),\\
                            \ea
\end{equation*}
where $U=U_1^\top U_2, V=V_1^\top V_2$ are clearly orthogonal
matrices.
Now let us derive an upper bound for $\Tr(Y_1^\top
Y_2-\bar{Y}_1^\top\bar{Y}_2)$. It is known that an orthogonal matrix
$U$ is a maximizing matrix for the problem
\begin{align*} \max\{\Tr(AU): U \mbox{ is orthogonal} \}\end{align*}
if and only if $AU$ is positive semidefinite matrix (see 7.4.9 in
\citep{Horn-Johnson-book-1985}). It is also known that when $AB$ is
positive semidefinite, \bea\label{eq:trace-singular-value-lemma-1}
\Tr(AB) = \sum_i \sigma_i(AB) \leq \sum_i\sigma_i(A)\sigma_i(B).\eea
Thus, $\Tr((\Sigma - \bar\Sigma)^\top U (\Gamma-\bar\Gamma)V^\top)$,
$\Tr((\Sigma-\bar\Sigma)^\top U\bar\Gamma V^\top)$ and
$\Tr(\bar\Sigma U(\Gamma-\bar\Gamma)V^\top)$ achieve their maximum,
if and only if $(\Sigma - \bar\Sigma)^\top U
(\Gamma-\bar\Gamma)V^\top$, $(\Sigma-\bar\Sigma)^\top U\bar\Gamma
V^\top$ and $\bar\Sigma U(\Gamma-\bar\Gamma)V^\top$ are all positive
semidefinite. Applying \eqref{eq:trace-singular-value-lemma-1} to
these three terms, we get $\Tr((\Sigma - \bar\Sigma)^\top U
(\Gamma-\bar\Gamma)V^\top)\leq\sum_i\sigma_i(\Sigma-\bar\Sigma)\sigma_i(\Gamma-\bar\Gamma)$,
$\Tr((\Sigma-\bar\Sigma)^\top U\bar\Gamma
V^\top)\leq\sum_i\sigma_i(\Sigma-\bar\Sigma)\sigma_i(\bar\Gamma)$
and $\Tr(\bar\Sigma
U(\Gamma-\bar\Gamma)V^\top)\leq\sum_i\sigma_i(\bar\Sigma)\sigma_i(\Gamma-\bar\Gamma).$
Thus, without loss of generality, assuming $k\leq l\leq s\leq t$, we
have,
\bean\ba{cl}&\|Y_1-Y_2\|_F^2-\|S_\nu(Y_1)-S_\nu(Y_2)\|_F^2\\
\geq &
\dis\sum_{i=1}^s\sigma_i^2-\sum_{i=1}^k(\sigma_i-\nu)^2+\sum_{i=1}^t\gamma_i^2-\sum_{i=1}^l(\gamma_i-\nu)^2-
2(\sum_{i=1}^l\sigma_i\nu+\sum_{i=l+1}^s\sigma_i\gamma_i+\sum_{i=1}^k(\gamma_i-\nu)\nu+\sum_{i=k+1}^l\sigma_i(\gamma_i-\nu))\\
=&\dis\sum_{i=k+1}^l(2\gamma_i\nu-\nu^2+\sigma_i^2-2\sigma_i\gamma_i)+(\sum_{i=l+1}^s\sigma_i^2+\sum_{i=l+1}^t\gamma_i^2-\sum_{i=l+1}^s2\sigma_i\gamma_i).\ea\eea

Now
$$\sum_{i=l+1}^s\sigma_i^2+\sum_{i=l+1}^t\gamma_i^2-\sum_{i=l+1}^s2\sigma_i\gamma_i\geq
0$$ since $t\geq s$ and $\sigma_i^2+\gamma_i^2-2\sigma_i\gamma_i\geq
0$. Also, since the function $f(x):=2\gamma_ix-x^2$ is monotonely
increasing in $(-\infty,\gamma_i]$ and
$\sigma_i<\nu\leq\gamma_i,i=k+1,\ldots,l$,
$$2\gamma_i\nu-\nu^2+\sigma_i^2-2\sigma_i\gamma_i> 0,i=k+1,\ldots,l.$$
Thus we get \bean
D(Y_1,Y_2):=\|Y_1-Y_2\|_F^2-\|S_\nu(Y_1)-S_\nu(Y_2)\|_F^2\geq 0;\eea
i.e., \eqref{ieq:decrese-shrinkage} holds.


Also, $D(Y_1,Y_2)$ achieves its minimum value if and only if
$\Tr((\Sigma - \bar\Sigma)^\top U (\Gamma-\bar\Gamma)V^\top)$,
$\Tr((\Sigma-\bar\Sigma)^\top U\bar\Gamma V^\top)$ and
$\Tr(\bar\Sigma U(\Gamma-\bar\Gamma)V^\top)$ achieve their maximum
values simultaneously.

Furthermore, if equality in \eqref{ieq:decrese-shrinkage} holds,
i.e., $D(Y_1,Y_2)$ achieves its minimum, and its minimum is zero,
then $k=l$, $s=t$, and $\sigma_i=\gamma_i,i=k+1,\ldots,s$, which
further implies $\Sigma-\bar\Sigma=\Gamma-\bar\Gamma$ and
$\Tr((\Sigma - \bar\Sigma)^\top U (\Gamma-\bar\Gamma)V^\top)$
achieves its maximum. By applying the result 7.4.13 in
\citep{Horn-Johnson-book-1985}, we get \bean \Sigma - \bar\Sigma = U
(\Gamma-\bar\Gamma)V^\top, \eea which further implies that
\begin{equation}\label{eq:last-eq-proof-lemma-1}
Y_1-Y_2=S_\nu(Y_1)-S_\nu(Y_2).\end{equation} 
To conclude, clearly $\|S_\nu(Y_1)-S_\nu(Y_2)\|_F=\|Y_1-Y_2\|_F$ if
\eqref{eq:last-eq-proof-lemma-1} holds. \qed
%

The following two lemmas and theorem and their proofs are analogous
to results and their proofs in Hale \etal
\citep{Hale-Yin-Zhang-2007}.
\begin{lemma}\label{lemma:h-non-expansive}
Let $\A X = A\vvec(X)$ and assume that
$\tau\in(0,2/\lambda_{max}(A^\top A))$. 
Then the operator $h(\cdot) = I(\cdot) -\tau g(\cdot)$ is
non-expansive, i.e., $\|h(X)-h(X')\|_F\leq\|X-X'\|_F$. Moreover,
$h(X)-h(X')=X-X'$ if and only if $\|h(X)-h(X')\|_F=\|X-X'\|_F$.
\end{lemma}
\proof First, we note that since $\tau\in(0,2/\lambda_{max}(A^\top
A))$, $-1<\lambda_i(I-\tau A^\top A)\leq 1,\forall i,$ where
$\lambda_i(I-\tau A^\top A)$ is the $i$-th eigenvalue of $I-\tau
A^\top A$. Hence,
\begin{align*}
\|h(X)-h(X')\|_F & = \|(I-\tau A^\top A)(\vvec(X)-\vvec(X'))\|_2
\leq \|I-\tau A^\top A\|_2 \|\vvec(X)-\vvec(X')\|_2 \\
& \leq \|\vvec(X)-\vvec(X')\|_2 = \|X-X'\|_F.
\end{align*}
Moreover, $\|h(X)-h(X')\|_F=\|X-X'\|_F$ if and only if the
inequalities above are equalities, which happens if and only if
\begin{align*}(I-\tau A^\top A)(\vvec(X)-\vvec(X')) = \vvec(X)-\vvec(X'),\end{align*}
i.e., if and only if $h(X)-h(X')=X-X'.$
\qed


\begin{lemma}\label{lemma:equal-distance-fixed-point}
Let $X^*$ be an optimal solution to problem
\eqref{prob:min-unconstrained}, $\tau\in(0,2/\lambda_{max}(A^\top
A))$ and $\nu=\tau\mu$. Then $X$ is also an optimal solution to
problem \eqref{prob:min-unconstrained} if and only if
\bea\label{eq:1st-equation-lemma-equal-distance-fixed-point}\|S_\nu(
h(X))-S_\nu(h(X^*))\|_F\equiv\|S_\nu(h(X))-X^*\|_F=\|X-X^*\|_F.\eea
\end{lemma}
\proof The ``only if'' part is an immediate consequence of Corollary
\ref{corollary:1}.  For the ``if'' part, from Lemmas
\ref{lemma:shrinkage non-expansive} and \ref{lemma:h-non-expansive},
\begin{align*}\|X-X^*\|_F=\|S_\nu(
h(X))-S_\nu(h(X^*))\|_F\leq\|h(X)-h(X^*)\|_F\leq\|X-X^*\|_F.\end{align*}
Hence, both inequalities hold with equality. Therefore, first using
Lemma \ref{lemma:shrinkage non-expansive} and then Lemma
\ref{lemma:h-non-expansive} we obtain
\begin{align*}S_\nu(h(X))-S_\nu( h(X^*))=h(X)-h(X^*)=X-X^*,\end{align*}
which implies $S_\nu(h(X))=X$ since $S_\nu( h(X^*))=X^*$. It then
follows from Corollary \ref{corollary:1} that $X$ is an optimal
solution to problem \eqref{prob:min-unconstrained}. \qed


We now claim that the fixed-point iterations
\eqref{eq:fpc-one-step-scheme} converge to an optimal solution of
problem \eqref{prob:min-unconstrained}.
\begin{theorem}\label{the:global-convergence}
The sequence $\{X^{k}\}$ generated by the fixed point iterations
with $\tau\in(0,2/\lambda_{max}(A^\top A))$ converges to some
$X^*\in\X^*,$ where $\X^*$ is the set of optimal solutions of
problem \eqref{prob:min-unconstrained}.
\end{theorem}

\proof Since both $S_\nu(\cdot)$ and $h(\cdot)$ are non-expansive,
$S_\nu(h(\cdot))$ is also non-expansive. Therefore, $\{X^k\}$ lies
in a compact set and must have a limit point, say
$\bar{X}=\lim_{j\rightarrow\infty}X^{k_j}.$ Also, for any
$X^*\in\X^*$,
\begin{align*}\|X^{k+1}-X^*\|_F=\|S_\nu(h(X^k))-S_\nu(h(X^*))\|_F\leq
\|h(X^k)-h(X^*)\|_F\leq\|X^k-X^*\|_F, \end{align*} which means that
the sequence $\{\|X^k-X^*\|_F\}$ is monotonically non-increasing.
Therefore,
\bea\label{eq:limit-distance-lemma-global-convergence}\lim_{k\rightarrow\infty}\|X^k-X^*\|_F=\|\bar{X}-X^*\|_F,\eea
where $\bar{X}$ can be any limit point of $\{X^k\}$. By the
continuity of $S_\nu(h(\cdot))$, the image of $\bar{X}$,
$$S_\nu(h(\bar{X}))=\lim_{j\rightarrow\infty} S_\nu(h(X^{k_j}))=\lim_{j\rightarrow\infty}X^{k_j+1},$$ is also a limit
point of $\{X^k\}$. Therefore, we have
$$\|S_\nu(h(\bar{X}))-S_\nu(h(X^*))\|_F=\|S_\nu(h(\bar{X}))-X^*\|_F=\|\bar{X}-X^*\|_F,$$
which allows us to apply Lemma
\ref{lemma:equal-distance-fixed-point} to get that $\bar{X}$ is an
optimal solution to problem \eqref{prob:min-unconstrained}.




Finally, by setting $X^*=\bar{X}\in\X^*$ in
\eqref{eq:limit-distance-lemma-global-convergence}, we get that
\begin{equation*}\lim_{k\rightarrow\infty}\|X^k-\bar{X}\|_F=\lim_{j\rightarrow\infty}\|X^{k_j}-\bar{X}\|_F=0,\end{equation*}
i.e., $\{X^k\}$ converges to its unique limit point $\bar{X}.$ \qed

\section{Fixed point continuation}
In this section, we discuss a continuation technique (i.e., homotopy
approach) for accelerating the convergence of the fixed point
iterative algorithm \eqref{eq:fpc-one-step-scheme}.
\subsection{Continuation}
Inspired by the work of Hale \etal \citep{Hale-Yin-Zhang-2007}, we
first describe a continuation technique to accelerate the
convergence of the fixed point iteration
\eqref{eq:fpc-one-step-scheme}. Our fixed point continuation (FPC)
iterative scheme for solving \eqref{prob:min-unconstrained} is
outlined below.
\begin{algorithm}\label{alg:fpc}{Fixed Point Continuation (FPC)}
\begin{itemize}
\item Initialize: Given $X_0$, $\bar\mu>0$. Select
$\mu_1>\mu_2>\cdots>\mu_L=\bar\mu>0.$ Set $X=X_0$.
\item {\bf for} $\mu=\mu_1,\mu_2,\ldots,\mu_L$, {\bf do}
\begin{itemize}
\item {\bf while} NOT converged, {\bf do}
\begin{itemize}
    \item select $\tau>0$
    \item compute $Y = X-\tau\A^*(\A(X)-b)$, and SVD of
    $Y$, $Y=U\Diag(\sigma) V^\top$
    \item compute $X=U\Diag(s_{\tau\mu}(\sigma))V^\top$ 
\end{itemize}
\item {\bf end while}
\end{itemize}
{\bf end for}
\end{itemize}
\end{algorithm}
The parameter $\eta_\mu$ determines the rate of reduction of the
consecutive $\mu_k$, i.e.,
$$\mu_{k+1}=\max\{\mu_k\eta_\mu,\bar\mu\},\qquad k=1,\ldots,L-1$$

\subsection{Stopping criteria for inner iterations} Note that in the fixed point continuation algorithm, in the $k$-th inner iteration we
solve problem \eqref{prob:min-unconstrained} for a fixed
$\mu=\mu_k$. There are several ways to determine when to stop this
inner iteration, decrease $\mu$ and go to the next inner iteration.
The optimality conditions for \eqref{prob:min-unconstrained} is
given by \eqref{optcond:exist-W-unconstrained-problem-a} and
\eqref{optcond:exist-W-unconstrained-problem-b}. Thus we can use the
following condition as a stopping criterion:
\bea\label{stopping-rule-g}\|U_kV_k^\top+g^k/\mu\|_2 - 1 < gtol,
\eea where $gtol$ is a small positive parameter. However, the
expense of computing the largest singular value of a large matrix
greatly decreases the speed of the algorithm. Hence, we do not use
this criterion as a stopping rule for large matrices. Instead, we
use the criterion
\bea\label{stopping-rule-x}\frac{\|X^{k+1}-X^k\|_F}{\max\{1,\|X^k\|_F\}}
< xtol,\eea where $xtol$ is a small positive number, since when
$X^k$ gets close to an optimal solution $X^*$, the distance between
$X^k$ and $X^{k+1}$ should become very small. 

\subsection{Debiasing}
Debiasing is another technique that can improve the performance of
FPC. Debiasing has been used in compressed sensing algorithms for
solving \eqref{prob:basis-pursuit-intro} and its variants, where
debiasing is performed after a support set $\mathcal{I}$ has been
tentatively identified. Debiasing is the process of solving a least
squares problem restricted to the support set $\mathcal{I}$, i.e.,
we solve \bea\label{prob:debias-CS}\min &
\|A_\mathcal{I}x_\mathcal{I}-b\|_2, \eea where $A_\mathcal{I}$ is a
submatrix of $A$ whose columns correspond to the support index set
$\mathcal{I}$, and $x_\mathcal{I}$ is a subvector of $x$
corresponding to $\mathcal{I}$.

Our debiasing procedure for the matrix completion problem differs
from the procedure used in compressed sensing since the concept of a
support set is not applicable. When we do debiasing, we fix the
matrices $U^k$ and $V^k$ in the singular value decomposition of
$X^k$ and then solve a least squares problem to determine the
correct singular values $\sigma\in\br_+^r$; i.e., we solve
\bea\label{prob:debias-matrix-completion}\min_{\sigma\geq 0} &
\|\A(U^k\Diag(\sigma)V^{k^\top})-b\|_2, \eea where $r$ is the rank
of current matrix $X^k$. Because debiasing can be costly, we use a
rule proposed in \citep{Wen-Yin-Goldfarb-Zhang-2009} to decide when
to do it. In the continuation framework, we know that in each
subproblem with a fixed $\mu$, $\|X_{k+1}-X_k\|_F$ converges to
zero, and $\|g\|_2$ converges to $\mu$ when $X_{k}$ converges to the
optimal solution of the subproblem. We therefore choose to do
debiasing when $\|g\|_2/\|X_{k+1}-X_k\|_F$ becomes large because
this indicates that the change between two consecutive iterates is
relatively small. Specifically, we call for debiasing in the solver
FPC3 (see Section 7) when $\|g\|_2/\|X_{k+1}-X_k\|_F
> 10.$

\section{Bregman iterative algorithm} Algorithm FPC is designed to
solve \eqref{prob:min-unconstrained}, an optimal solution of which
approaches an optimal solution of the nuclear norm minimization
problem \eqref{prob:nuclear-norm-minimization-intro} as $\mu$ goes
to zero. However, by incorporating FPC into a Bregman iterative
technique, we can solve \eqref{prob:nuclear-norm-minimization-intro}
by solving a limited number of instances of
\eqref{prob:min-unconstrained}, each corresponding to a different
$b$.

Given a convex function $J(\cdot)$, the Bregman distance
\citep{Bregman-1967} of the point $u$ from the point $v$ is defined
as
\bea\label{definition:bregman-distance}D_J^p(u,v):=J(u)-J(v)-<p,u-v>,\eea
where $p\in\partial J(v)$ is some subgradient in the subdifferential
of $J$ at the point $v$.

Bregman iterative regularization was introduced by Osher \etal in
the context of image processing
\citep{Osher-Burger-Goldfarb-Xu-Yin-2005}. Specifically, in
\citep{Osher-Burger-Goldfarb-Xu-Yin-2005}, the Rudin-Osher-Fatemi
\citep{Rudin-Osher-Fatemi-1992} model
\bea\label{prob:ROF}u=\argmin_u\mu\int |\nabla
u|+\df{1}{2}\|u-b\|_2^2\eea was extended to an iterative
regularization model by replacing the total variation functional
$$J(u)=\mu TV(u)=\mu\int |\nabla u|,$$ by the Bregman distance with respect to $J(u)$.
This Bregman iterative regularization procedure recursively solves
\bea\label{iterative-Osher-Bregman-u_k+1}
u^{k+1}\leftarrow\min_uD_J^{p^k}(u,u^k)+\df12\|u-b\|_2^2\eea for
$k=0,1,\ldots$ starting with $u^0=\zero$ and $p^0=\zero$. Since
\eqref{iterative-Osher-Bregman-u_k+1} is a convex programming
problem, the optimality conditions are given by $\zero\in\partial
J(u^{k+1})-p^k+u^{k+1}-b,$ from which we get the update formula for
$p^{k+1}:$
\bea\label{eq:update-p_k+1-Bregman}p^{k+1}:=p^k+b-u^{k+1}.\eea
Therefore, the Bregman iterative scheme is given by
\bea\label{iterative-ROF-Bregman-MMS-version1}\left\{\ba{l}u^{k+1}\leftarrow\min_uD_J^{p^k}(u,u^k)+\df12\|u-b\|_2^2\\p^{k+1}=p^k+b-u^{k+1}.\ea\right.\eea
Interestingly, this turns out to be equivalent to the iterative
process
\bea\label{iterative-ROF-Bregman-MMS-version2}\left\{\ba{l}b^{k+1}=b+(b^k-u^k)\\u^{k+1}\leftarrow\min_{u}J(u)+\df{1}{2}\|u-b^{k+1}\|_2^2,\ea\right.\eea
which can be easily implemented using existing algorithms for
\eqref{prob:ROF} with different inputs $b$.

Subsequently, Yin \etal \citep{Yin-Osher-Goldfarb-Darbon-2008}
proposed solving the basis pursuit problem
\eqref{prob:basis-pursuit-intro} by applying the Bregman iterative
regularization algorithm to
\bea\label{prob:basis-pursuit-unconstrained} \min_x
J(x)+\df12\|Ax-b\|_2^2\eea for $J(x)=\mu\|x\|_1,$ and obtained the
following two equivalent iterative schemes analogous to
\eqref{iterative-ROF-Bregman-MMS-version1} and
\eqref{iterative-ROF-Bregman-MMS-version2}, respectively:
\begin{itemize}
\item Version 1:
\begin{itemize}
\item $x^0\leftarrow \zero, p^0\leftarrow \zero,$
\item for $k=0,1,\ldots$ do
\item $x^{k+1}\leftarrow\argmin_x D_J^{p^k}(x,x^k)+\df12\|Ax-b\|_2^2$
\item $p^{k+1}\leftarrow p^k-A^\top(Ax^{k+1}-b)$
\end{itemize}

\item Version 2:
\begin{itemize}
\item $b^0\leftarrow\zero, x^0\leftarrow \zero,$
\item for $k=0,1,\ldots$ do
\item $b^{k+1}\leftarrow b+(b^k-Ax^k)$
\item $x^{k+1}\leftarrow\argmin_xJ(x)+\df12\|Ax-b^{k+1}\|_2^2.$
\end{itemize}
\end{itemize}

One can also use the Bregman iterative regularization algorithm
applied to the unconstrained problem \eqref{prob:min-unconstrained}
to solve the nuclear norm minimization problem
\eqref{prob:nuclear-norm-minimization-intro}. That is, one
iteratively solves \eqref{prob:min-unconstrained} by
\bea\label{eq:bregman-u_k+1-nuclear-norm}X^{k+1}\leftarrow\min_X
D_J^{p^k}(X,X^k)+\df12\|\A(X)-b\|_2^2,\eea and updates the
subgradient $p^{k+1}$ by
\bea\label{eq:bregman-p_k+1-update-nuclear-norm}p^{k+1}:=p^k-\A^*(A(X^{k+1})-b),\eea
where $J(X)=\mu\|X\|_*$.

Equivalently, one can also use the following iterative scheme:
\bea\label{iterative-Bregman-2-nuclear-norm}\left\{\ba{l}b^{k+1}\leftarrow
b+(b^k-\A(X^k))\\X^{k+1}\leftarrow\arg\min_X\mu\|X\|_*+\df{1}{2}\|\A(X)-b^{k+1}\|_2^2.\ea\right.\eea

Thus, our Bregman iterative algorithm for nuclear norm minimization
\eqref{prob:nuclear-norm-minimization-intro} can be outlined as
follows.
\begin{algorithm}\label{alg:bregman}{Bregman Iterative Algorithm}
\begin{itemize}
\item $b^0\leftarrow\zero, X^0\leftarrow\zero,$
\item for $k=0,1,\ldots$ do
\item $b^{k+1}\leftarrow b+(b^k-\A(X^k)),$
\item
$X^{k+1}\leftarrow\arg\min_X\mu\|X\|_*+\df{1}{2}\|\A(X)-b^{k+1}\|_2^2$.
\end{itemize}
\end{algorithm}
The last step can be solved by Algorithm FPC.

\section{An approximate SVD based FPC algorithm: FPCA}\label{sec:fast-svd} Computing singular value decompositions is the main
computational cost in Algorithm FPC. Consequently, instead of
computing the full SVD of the matrix $Y$ in each iteration, we
implemented a variant of algorithm FPC in which we compute only a
rank-$r$ approximation to $Y$, where $r$ is a predetermined
parameter. We call this approximate SVD based FPC algorithm (FPCA).
This approach greatly reduces the computational effort required by
the algorithm. Specifically, we compute an approximate SVD by a fast
Monte Carlo algorithm: the Linear Time SVD algorithm developed by
Drineas \etal \citep{Drineas-Kannan-Mahoney-2006}. For a given
matrix $A\in\mathbb{R}^{m\times n}$, and parameters
$c_s,k_s\in\mathbb{Z}^+$ with $1\leq k_s\leq c_s\leq n$ and
$\{p_i\}_{i=1}^n,$ $p_i\geq 0, \sum_{i=1}^np_i=1$, this algorithm
returns an approximation to the largest $k_s$ singular values and
corresponding left singular vectors of the matrix $A$ in linear
$O(m+n)$ time. The Linear Time SVD Algorithm is outlined below.
\begin{algorithm}\label{alg:Linear-Time-SVD}{Linear Time Approximate SVD
Algorithm}\citep{Drineas-Kannan-Mahoney-2006}
\begin{itemize}
\item {\bf Input:} $A\in\mathbb{R}^{m\times n}$, $c_s,k_s\in\mathbb{Z}^+$ \st $1\leq k_s\leq c_s\leq
n$, $\{p_i\}_{i=1}^n$ \st $p_i\geq 0, \sum_{i=1}^np_i=1$.
\item {\bf Output:} $H_k\in\mathbb{R}^{m\times k_s}$ and $\sigma_t(C), t=1,\ldots,k_s.$
\begin{itemize}
\item For $t=1$ to $c_s$,
\begin{itemize}
\item Pick $i_t\in 1,\ldots,n$ with $Pr[i_t=\alpha] =
p_\alpha,\alpha=1,\ldots,n.$
\item Set $C^{(t)}=A^{(i_t)}/\sqrt{c_sp_{i_t}}.$
\end{itemize}
\item Compute $C^\top C$ and its SVD; say $
C^\top C=\sum_{t=1}^{c_s}\sigma_t^2(C)y^{t}{y^t}^\top.$
\item Compute $h^t=Cy^t/\sigma_t(C)$ for $t=1,\ldots,k_s.$
\item Return $H_{k_s}$, where $H_{k_s}^{(t)}=h^t,$ and $\sigma_t(C),
t=1,\ldots,k_s.$
\end{itemize}
\end{itemize}
\end{algorithm}

The outputs $\sigma_t(C),t=1,\ldots,k_s$ are approximations to the
largest $k_s$ singular values and $H_{k_s}^{(t)},t=1,\ldots,k$ are
approximations to the corresponding left singular vectors of the
matrix $A$. Thus, the SVD of $A$ is approximated by $$A\approx
A_{k_s}:= H_{k_s}\Diag(\sigma(C))(A^\top
H_{k_s}\Diag(1/\sigma(C))^\top.$$ Drineas
\etal\citep{Drineas-Kannan-Mahoney-2006} prove that with high
probability, the following estimate holds for both $\xi=2$ and
$\xi=F$:
\bea\label{eq:error-fast-svd}\|A-A_{k_s}\|_\xi^2\leq\min_{D:\rank(D)\leq
k_s}\|A-D\|_\xi^2+poly(k_s,1/c_s)\|A\|_F^2,\eea where
$poly(k_s,1/c_s)$ is a polynomial in $k_s$ and $1/c_s$. Thus,
$A_{k_s}$ is a approximation to the best rank-$k_s$ approximation to
$A$. (For any matrix $M\in\br^\mtn$ with SVD
$M=\sum_{i=1}^r\sigma_iu_iv_i^\top$, where
$\sigma_1\geq\ldots\geq\sigma_r>0,u_i\in\br^m,v_i\in\br^n$, the best
rank-k approximation to $M$ is given by
$\bar{M}=\sum_{i=1}^k\sigma_iu_iv_i^\top$).

Note that in this algorithm, we compute an exact SVD of a smaller
matrix $C^\top C\in\mathbb{R}^{c_s\times c_s}$. Thus, $c_s$
determines the speed of this algorithm. If we choose a large $c_s$,
we need more time to compute the SVD of $C^\top C$. However, the
larger $c_s$ is, the more likely are the
$\sigma_t(C),t=1,\ldots,k_s$ to be close to the largest $k_s$
singular values of the matrix $A$ since the second term in the right
hand side of \eqref{eq:error-fast-svd} is smaller. In our numerical
experiments, we found that we could choose a relatively small $c_s$
so that the computational time was reduced without significantly
degrading the accuracy. In our tests, we obtained very good results
by choosing $c_s=2r_m-2$, where $r_m=\lfloor
(m+n-\sqrt{(m+n)^2-4p})/2\rfloor$ is, for a given number of entries
sampled, the largest rank of $\mtn$ matrices for which the matrix
completion problem has a unique solution.

There are many ways to choose the probabilities $p_i$. In our
numerical experiments in Section 7, we used the simplest one, i.e.,
we set all $p_i$ equal to $1/n$. For other choices of $p_i$, see
\citep{Drineas-Kannan-Mahoney-2006} and the references therein.

In our numerical experiments, we set $k_s$ using the following
procedure. In the $k$-th iteration, when computing the approximate
SVD of $Y^k=X^k-\tau g^k$, we set $k_s$ equal to the number of
components in $\bar{s}_{k-1}$ that are no less than
$\epsilon_{k_s}\max\{\bar{s}_{k-1}\},$ where $\epsilon_{k_s}$ is a
small positive number and $\max\{\bar{s}_{k-1}\}$ is the largest
component in the vector $\bar{s}_{k-1}$ used to form
$X^k=U^{k-1}\Diag(\bar{s}_{k-1}){V^{k-1}}^\top$. Note that $k_s$ is
non-increasing in this procedure. However, if $k_s$ is too small at
some iteration, the non-expansive property
\eqref{ieq:decrese-shrinkage} of the shrinkage operator $S_\nu$ may
be violated since the approximate SVD is not a valid approximation
when $k_s$ is too small. Thus, in algorithm FPCA, if
\eqref{ieq:decrese-shrinkage} is violated 10 times, we increase
$k_s$ by $1$. Our numerical experience indicates that this technique
makes our algorithm very robust.

Our numerical results in Section 7 show that this approximate SVD
based FPC algorithm: FPCA, is very fast, robust, and significantly
outperforms other solvers (such as SDPT3) in recovering low-rank
matrices. This result is not surprising. One reason for this is that
in the approximate SVD algorithm, we compute a low-rank
approximation to the original matrix. Hence, the iterative matrices
produced by our algorithm are more likely to be of low-rank than an
exact solution to the nuclear norm minimization problem
\eqref{prob:nuclear-norm-matrix-completion-intro}, or equivalently,
to the SDP \eqref{prob:sdp-nuclear-norm-primal}, which is exactly
what we want. Some convergence/recoverability properties of a
variant of FPCA, which uses a truncated SVD rather than a randomized
SVD at each step, are discussed in \citep{Goldfarb-Ma-2009}.

\section{Numerical results}

In this section, we report on the application of our FPC, FPCA and
Bregman iterative algorithms to a series of matrix completion
problems of the form \eqref{prob:rank-matrix-completion-intro} to
demonstrate the ability of these algorithms to efficiently recover
low-rank matrices.

To illustrate the performance of our algorithmic approach combined
with exact and approximate SVD algorithms, different stopping rules,
and with or without debiasing, we tested the following solvers.
\begin{itemize}
\item FPC1. Exact SVD, no debiasing, stopping rule:
\eqref{stopping-rule-x}.
\item FPC2. Exact SVD, no debiasing, stopping rule:
\eqref{stopping-rule-g} and \eqref{stopping-rule-x}.
\item FPC3. Exact SVD with debiasing, stopping rule: \eqref{stopping-rule-x}.
\item FPCA. Approximate SVD, no debiasing, stopping rule: \eqref{stopping-rule-x}.
\item Bregman. Bregman iterative method using FPC2 to solve the
subproblems.
\end{itemize}

\subsection{FPC and Bregman iterative algorithms for random matrices}
In our first series of tests, we created random matrices
$M\in\br^{m\times n}$ with rank $r$ by the following procedure: we
first generated random matrices $M_L\in\br^{m\times r}$ and
$M_R\in\br^{n\times r}$ with i.i.d. Gaussian entries and then set
$M=M_LM_R^\top$. We then sampled a subset $\Omega$ of $p$ entries
uniformly at random. For each problem with $m\times n$ matrix $M$,
measurement number $p$ and rank $r$, we solved 50 randomly created
matrix completion problems. We use $SR=p/(mn)$, i.e., the number of
measurements divided by the number of entries of the matrix, to
denote the sampling ratio. We also list $FR=r(m+n-r)/p$, i.e. the
dimension of the set of rank $r$ matrices divided by the number of
measurements, in the tables. Note that if $FR>1$, then there is
always an infinite number of matrices with rank $r$ with the given
entries, so we cannot hope to recover the matrix in this situation.
We use $r_m$ to denote the largest rank such that $FR\leq 1$, i.e.,
$r_m=\lfloor (m+n-\sqrt{(m+n)^2-4p})/2\rfloor$. We use $NS$ to
denote the number of matrices that are recovered successfully. We
use $AT$ to denote the average time (seconds) for the examples that
are successfully solved.

We used the relative error \bean rel.err. :=
\frac{\|X_{opt}-M\|_F}{\|M\|_F}\eea to estimate the closeness of
$X_{opt}$ to $M$, where $X_{opt}$ is the ``optimal'' solution to
\eqref{prob:nuclear-norm-matrix-completion-intro} produced by our
algorithms. We declared $M$ to be recovered if the relative error
was less than $10^{-3}$, which is the criterion used in
\citep{Recht-Fazel-Parrilo-2007} and \citep{Candes-Recht-2008}. We
use $RA, RU, RL$ to denote the average, largest and smallest
relative error of the successfully recovered matrices, respectively.

We summarize the parameter settings used by the algorithms in Table
\ref{table:parameters-fpc}. We use $I_m$ to denote the maximum
number of iterations allowed for solving each subproblem in FPC,
i.e., if the stopping rules \eqref{stopping-rule-x} (and
\eqref{stopping-rule-g}) are not satisfied after $I_m$ iterations,
we terminate the subproblem and decrease $\mu$ to start the next
subproblem.
\begin{table}[ht]
\begin{center}\caption{Parameters in
Algorithm FPC}\label{table:parameters-fpc}
\begin{tabular}{c|c}\hline\hline
FPC &
{$\bar{\mu}=10^{-8},\eta_{\mu}=1/4,\mu_1=\eta_{\mu}\|\A^*b\|_2,\tau=1,xtol=10^{-10},gtol=10^{-4},I_m=500,X_0=\mathbf{0}$}
\\\hline
Approx SVD & $c_s=2r_m-2,\epsilon_{k_s}=10^{-2},p_i=1/n,\forall i$
\\\hline
\end{tabular}
\end{center}
\end{table}

All numerical experiments were run in MATLAB 7.3.0 on a Dell
Precision 670 workstation with an Intel Xeon(TM) 3.4GHZ CPU and 6GB
of RAM.

The comparisons between FPC1, FPC2, FPC3 and SDPT3 for small matrix
completion problems are presented in Table
\ref{table:Num-res-fpcE-sdpt3-random-small}. From Table
\ref{table:Num-res-fpcE-sdpt3-random-small} we can see that FPC1 and
FPC2 achieve almost the same recoverability and relative error,
which means that as long as we set $xtol$ to be very small (like
$10^{-10}$ ), we only need to use \eqref{stopping-rule-x} as the
stopping rule for the inner iterations. That is, use of stopping
rule \eqref{stopping-rule-g} does not affect the performance of the
algorithm. Of course FPC2 costs more time than FPC1 since more
iterations are sometimes needed to satisfy the stopping rules in
FPC2. While FPC3 can improve the recoverability, it costs more time
for performing debiasing. SDPT3 seems to obtain more accurate
solutions than FPC1, FPC2 or FPC3.


\begin{table}[ht]
\begin{center}\caption{Comparisons of FPC1, FPC2, FPC3 and SDPT3 for randomly created small matrix
completion problems (m=n=40, p=800,
SR=0.5)}\label{table:Num-res-fpcE-sdpt3-random-small}
\begin{tabular}{ c | c | c | c c c c c }\hline
r &  FR  & Solver & NS & AT & RA & RU & RL  \\\hline\hline

1 &  0.0988 & FPC1 & 50 & 1.81 & 1.67e-9 & 1.22e-8 & 6.06e-10 \\
  &         & FPC2 & 50 & 3.61 & 1.32e-9 & 1.20e-8 & 2.55e-10 \\
  &         & FPC3 & 50 & 16.81 & 1.06e-9 & 2.22e-9 & 5.68e-10 \\
  &         & SDPT3 & 50 & 1.81 & 6.30e-10 & 3.46e-9 & 8.72e-11 \\\hline

2 & 0.1950 & FPC1 & 42 & 3.05 & 1.01e-6 & 4.23e-5 & 8.36e-10 \\
  &        & FPC2 & 42 & 17.97 & 1.01e-6 & 4.23e-5 & 2.78e-10 \\
  &        & FPC3 & 49 & 16.86 & 1.26e-5 & 3.53e-4 & 7.62e-10 \\
  &        & SDPT3 & 44 & 1.90 & 1.50e-9 & 7.18e-9 & 1.82e-10 \\\hline

3 & 0.2888 & FPC1 & 35 & 5.50 & 9.72e-9 & 2.85e-8 & 1.93e-9 \\
  &        & FPC2 & 35 & 20.33 & 2.17e-9 & 1.41e-8 & 3.88e-10 \\
  &        & FPC3 & 42 & 16.87 & 3.58e-5 & 7.40e-4 & 1.34e-9 \\
  &        & SDPT3 & 37 & 1.95 & 2.66e-9 & 1.58e-8 & 3.08e-10 \\\hline

4 & 0.3800 & FPC1 & 22 & 9.08 & 7.91e-5 & 5.46e-4 & 3.57e-9 \\
  &        & FPC2 & 22 & 18.43 & 7.91e-5 & 5.46e-4 & 4.87e-10 \\
  &        & FPC3 & 29 & 16.95 & 3.83e-5 & 6.18e-4 & 2.57e-9 \\
  &        & SDPT3 & 29 & 2.09 & 1.18e-8 & 7.03e-8 & 7.97e-10 \\\hline

5 & 0.4688 & FPC1 & 1 & 10.41 & 2.10e-8 & 2.10e-8 & 2.10e-8 \\
  &        & FPC2 & 1 & 17.88 & 2.70e-9 & 2.70e-9 & 2.70e-9 \\
  &        & FPC3 & 5 & 16.70 & 1.78e-4 & 6.73e-4 & 6.33e-9 \\
  &        & SDPT3 & 8 & 2.26 & 1.83e-7 & 8.12e-7 & 2.56e-9 \\\hline

6 & 0.5550 & FPC1 & 0 & --- & --- & --- & ---  \\
  &        & FPC2 & 0 & --- & --- & --- & --- \\
  &        & FPC3 & 0 & --- & --- & --- & --- \\
  &        & SDPT3 & 1 & 2.87 & 6.58e-7 & 6.58e-7 & 6.58e-7 \\\hline
\end{tabular}
\end{center}
\end{table}

To illustrate the performance of our Bregman iterative algorithm, we
compare the results of using it versus using FPC2 in Table
\ref{table:Numerical-results-matrix-completion-Bregman}. From our
numerical experience, for those problems for which the Bregman
iterative algorithm greatly improves the recoverability, the Bregman
iterative algorithm usually takes 2 to 3 iterations. Thus, in our
numerical tests, we fixed the number of subproblems solved by our
Bregman algorithm to 3. Since our Bregman algorithm achieves as good
a relative error as the FPC algorithm, we only report how many of
the examples that are successfully recovered by FPC, are improved
greatly by using our Bregman iterative algorithm. In Table
\ref{table:Numerical-results-matrix-completion-Bregman}, NIM is the
number of examples that the Bregman iterative algorithm outperformed
FPC2 greatly (the relative errors obtained from FPC2 were at least
$10^4$ times larger than those obtained by the Bregman algorithm).
From Table \ref{table:Numerical-results-matrix-completion-Bregman}
we can see that for more than half of the examples successfully
recovered by FPC2, the Bregman iterative algorithm improved the
relative errors greatly (from [$10^{-10}$, $10^{-9}$] to
[$10^{-16}$, $10^{-15}$]). Of course the run times for the Bregman
iterative algorithm were about three times that for algorithm FPC2,
since the former calls the latter three times to solve the
subproblems.

\begin{table}[ht]
\begin{center}\caption{Numerical results for the Bregman iterative method for small matrix completion
problems (m=n=40, p=800,
SR=0.5)}\label{table:Numerical-results-matrix-completion-Bregman}
\begin{tabular}{c c | c | c c | c c}\hline
\multicolumn{2}{c|}{Problem} & &\multicolumn{2}{c|}{FPC2} &
\multicolumn{2}{c}{Bregman}
\\\hline

r & FR  & NIM (NS)& RU & RL & RU & RL\\\hline\hline
1 & 0.0988 & 32 (50) & 2.22e-9 & 2.55e-10 & 1.87e-15 & 3.35e-16 \\
2 & 0.1950 & 29 (42) & 5.01e-9 & 2.80e-10 & 2.96e-15 & 6.83e-16 \\
3 & 0.2888 & 24 (35) & 2.77e-9 & 3.88e-10 & 2.93e-15 & 1.00e-15 \\
4 & 0.3800 & 10 (22) & 5.51e-9 & 4.87e-10 & 3.11e-15 & 1.30e-15
\\\hline
\end{tabular}
\end{center}
\end{table}

In the following, we discuss the numerical results obtained by our
approximate SVD based FPC algorithm (FPCA). We will see from these
numerical results that FPCA achieves much better recoverability and
is much faster than any of the solvers FPC1, FPC2, FPC3 or SDPT3.

We present the numerical results of FPCA for small (m=n=40) and
medium (m=n=100) problems in Tables
\ref{table:Num-res-fpc-sdp-matrix-completion-small},
and \ref{table:Num-res-fpc-sdp-matrix-completion-medium} 
respectively. Since we found that $xtol=10^{-6}$ is small enough to
guarantee very good recoverability, we set $xtol=10^{-6}$ in
algorithm FPCA and used only \eqref{stopping-rule-x} as stopping
rule for the inner iterations. From these tables, we can see that
our FPCA algorithm is much more powerful than SDPT3 for randomly
created matrix completion problems. When $m=n=40$ and $p=800$, and
the rank $r$ was less than or equal to 8, FPCA recovered the
matrices in all 50 examples. When rank $r=9$, it failed on only one
example. Even for rank $r=10$, which is almost the largest rank that
satisfies $FR\leq 1$, FPCA still recovered the solution in more than
$60\%$ of the examples. However, SDPT3 started to fail to recover
the matrices when the rank $r=2$. When $r=6$, there was only one
example out of 50 where the correct solution matrix was recovered.
When $r\geq 7$, none of the 50 examples could be recovered. For the
medium sized matrices $(m=n=100)$ we used $p=2000$, which is only a
$20\%$ measurement rate, FPCA recovered the matrices in all 50
examples when $r\leq 6$. For $r=7,$ FPCA recovered the matrices in
most of the examples (49 out of 50). When $r=8,$ more than $60\%$ of
the matrices were recovered successfully by FPCA. Even when $r=9,$
FPCA still recovered 1 matrices. However, SDPT3 could not recover
all of the matrices even when the rank $r=1$ and none of the
matrices were recovered when $r\geq 4.$ When we increased the number
of measurements to $3000$, we recovered the matrices in all 50
examples up to rank $r=12.$ When $r=13,14,$ we still recovered most
of them. However, SDPT3 started to fail for some matrices when
$r=3.$ When $r\geq 8$, SDPT3 failed to recover any of the matrices.
We can also see that for the medium sized problems, FPCA was much
faster than SDPT3.

\begin{table}[ht]
\begin{center}\caption{Numerical results for FPCA and SDPT3 for randomly created small matrix completion
problems (m=n=40, p=800,
SR=0.5)}\label{table:Num-res-fpc-sdp-matrix-completion-small}
\begin{tabular}{ c | c | c c c c c | c c c c c }\hline
\multicolumn{2}{c|}{Problems} & \multicolumn{5}{c|}{FPCA} &
\multicolumn{5}{c}{SDPT3}
\\\hline

r &  FR  & NS & AT & RA & RU & RL   & NS & AT & RA & RU & RL
\\\hline\hline
1 &  0.0988 & 50 & 3.49 & 3.92e-7 & 1.43e-6 & 2.72e-7 & 50 & 1.84 &
6.30e-10 & 3.46e-9 & 8.70e-11
\\

2 & 0.1950 & 50 & 3.60 & 1.44e-6 & 7.16e-6 & 4.41e-7 & 44 & 1.93 & 1.50e-9 & 7.18e-9 & 1.82e-10 \\

3 & 0.2888 & 50 & 3.97 & 1.91e-6 & 4.07e-6 & 9.28e-7 & 37 & 1.99
& 2.66e-9 & 1.58e-8 & 3.10e-10 \\

4 & 0.3800 & 50 & 4.03 & 2.64e-6 & 8.14e-6 & 1.54e-6& 29 & 2.12
& 1.18e-8 & 7.03e-8 & 8.00e-10 \\

5 & 0.4688 & 50 & 4.16 & 3.40e-6 & 7.62e-6 & 1.52e-6 & 8 & 2.30
& 1.83e-7 & 8.12e-7 & 2.60e-9 \\

6 & 0.5550 & 50 & 4.45 & 4.08e-6 & 7.62e-6 & 2.26e-6 & 1  & 2.89 & 6.58e-7 & 6.58e-7 & 6.58e-7 \\

7 & 0.6388 & 50 & 4.78 &
6.04e-6 & 1.57e-5 & 2.52e-6 & 0  & --- & --- & --- & --- \\

8 & 0.7200 & 50 & 4.99 & 8.48e-6 & 5.72e-5 & 3.72e-6 & 0  & --- & --- & --- & --- \\

9 & 0.7987 & 49 & 5.73 & 2.58e-5 & 5.94e-4 & 5.94e-6 & 0  & --- & --- & --- & --- \\

10 & 0.8750 & 30 & 7.20 & 8.64e-5 & 6.04e-4 & 8.48e-6 &  0  & --- & --- & --- & --- \\

11 & 0.9487 & 0  & --- & --- & --- & --- &  0  & --- & --- &
--- &
---
\\\hline
\end{tabular}
\end{center}
\end{table}

\begin{table}[ht]
\begin{center}\caption{Numerical results for FPCA and SDPT3 for randomly created medium matrix completion
problems
(m=n=100)}\label{table:Num-res-fpc-sdp-matrix-completion-medium}
\begin{tabular}{ c c c c | c c c c c | c c c c c}\hline
\multicolumn{4}{c|}{Problems} & \multicolumn{5}{c|}{FPCA} &
\multicolumn{5}{c}{SDPT3}
\\\hline

p & r & SR & FR & NS & AT & RA & RU & RL  & NS & AT & RA & RU & RL
\\\hline\hline
2000 & 1 & 0.2 & 0.0995 & 50 & 4.93 &
5.80e-6 & 1.53e-5 & 2.86e-6 & 47 & 15.10 & 1.55e-9 & 1.83e-8 & 1.40e-10 \\

2000 & 2 & 0.2 & 0.1980 & 50 & 5.26 &
6.10e-6 & 9.36e-6 & 4.06e-6 & 31 & 16.02 & 7.95e-9 & 8.69e-8 & 5.20e-10 \\

2000 & 3 & 0.2 & 0.2955 & 50 & 5.80 &
7.48e-6 & 1.70e-5 & 4.75e-6 & 13 & 19.23 & 1.05e-4 & 9.70e-4 & 9.08e-10 \\

2000 & 4 & 0.2 & 0.3920 & 50 & 9.33 & 1.09e-5 & 5.14e-5 & 6.79e-6 & 0  & --- & --- & --- & --- \\

2000 & 5 & 0.2 & 0.4875 & 50 & 5.42 & 1.61e-5 & 8.95e-5 & 8.12e-6 & 0  & --- & --- & --- & --- \\

2000 & 6 & 0.2 & 0.5820 & 50 & 7.02 & 2.62e-5 & 7.07e-5 & 8.72e-6 & 0  & --- & --- & --- & --- \\

2000 & 7 & 0.2 & 0.6755 & 49 & 8.69 & 7.69e-5 & 5.53e-4 & 1.11e-5 & 0  & --- & ---  & --- & --- \\

2000 & 8 & 0.2 & 0.7680 & 32 & 10.94 & 1.97e-4 & 8.15e-4 & 2.29e-5 & 0  & --- & --- & --- & --- \\

2000 & 9 & 0.2 & 0.8595 & 1  & 11.75 & 4.38e-4 & 4.38e-4 & 4.38e-4 & 0  & ---  & --- & --- & --- \\

2000 & 10 & 0.2 & 0.9500 & 0 & --- & --- & --- & --- & 0  & --- &
--- & --- & ---
\\\hline

3000 & 1 & 0.3 & 0.0663 & 50 & 7.73 & 1.97e-6 & 3.15e-6 & 1.22e-6 & 50 & 36.68 & 2.01e-10 & 9.64e-10 & 7.52e-11 \\

3000 & 2 & 0.3 & 0.1320 & 50 & 7.85 & 2.68e-6 & 8.41e-6 & 1.44e-6 & 50 & 36.50 & 1.13e-9 & 2.97e-9 & 1.77e-10 \\

3000 & 3 & 0.3 & 0.1970 & 50 & 8.10 & 2.82e-6 & 4.38e-6 & 1.83e-6 & 46 & 38.50 & 1.28e-5 & 5.89e-4 & 2.10e-10\\

3000 & 4 & 0.3 & 0.2613 & 50 & 8.94 & 3.57e-6 & 5.62e-6 & 2.64e-6 & 42 & 41.28 & 4.60e-6 & 1.21e-4 & 4.53e-10 \\

3000 & 5 & 0.3 & 0.3250 & 50 & 9.12 & 4.06e-6 & 8.41e-6 & 2.78e-6 & 32 & 43.92 & 7.82e-8 & 1.50e-6 & 1.23e-9 \\

3000 & 6 & 0.3 & 0.3880 & 50 & 9.24 & 4.84e-6 & 9.14e-6 & 3.71e-6 & 17 & 49.60 & 3.44e-7 & 4.29e-6 & 3.68e-9 \\

3000 & 7 & 0.3 & 0.4503 & 50 & 9.41 & 5.72e-6 & 1.09e-5 & 3.96e-6 & 3 & 59.18 & 1.43e-4 & 4.28e-4 & 1.57e-7 \\

3000 & 8 & 0.3 & 0.5120 & 50 & 9.62 & 6.37e-6 & 1.90e-5 & 4.43e-6 & 0 & --- & --- & --- & --- \\

3000 & 9 & 0.3 & 0.5730 & 50 & 10.35 & 6.32e-6 & 1.60e-5 & 4.56e-6 & 0 & --- & --- & --- & --- \\

3000 & 10 & 0.3 & 0.6333 & 50 & 10.93 & 8.45e-6 & 3.79e-5 & 5.59e-6 & 0 & --- & --- & --- & --- \\

3000 & 11 & 0.3 & 0.6930 & 50 & 11.58 & 1.41e-5 & 6.84e-5 & 6.99e-6 & 0 & --- & --- & --- & --- \\

3000 & 12 & 0.3 & 0.7520 & 50 & 12.17 & 1.84e-5 & 1.46e-4 & 8.84e-6 & 0 & --- & --- & --- & --- \\

3000 & 13 & 0.3 & 0.8103 & 48 & 15.24 & 5.12e-5 & 6.91e-4 & 1.25e-5 & 0 & --- & --- & --- & --- \\

3000 & 14 & 0.3 & 0.8680 & 39 & 18.85 & 2.35e-4 & 9.92e-4 & 2.05e-5 & 0 & --- & --- & --- & --- \\

3000 & 15 & 0.3 & 0.9250 & 0  & --- & --- & --- & --- & 0 & --- & --- & --- & --- \\

3000 & 16 & 0.3 & 0.9813 & 0  & --- & --- & --- & --- & 0 & --- &
--- & --- & ---
\\\hline
\end{tabular}
\end{center}
\end{table}

\subsection{Comparison of FPCA and SVT}

In this subsection we compare our FPCA algorithm against the SVT
algorithm proposed in \citep{Cai-Candes-Shen-2008}. The SVT code is
downloaded from {\em http://svt.caltech.edu}. We constructed two
sets of test problems. One set contained ``easy'' problems. These
problems are ``easy'' because the matrices are of very low-rank
compared to the matrix size and the number of samples, and hence
they are easy to recover. For all problems in this set, $FR$ was
less than 0.34. The other set contained ``hard'' problems, i.e.,
problems that are very challenging. These problems involved matrices
that are not of very low rank and for which sampled a very limited
number of entries. For this set of problems, $FR$ ranged from 0.40
to 0.87. The maximum iteration number in SVT was set to be 1000. All
other parameters were set to their default values in SVT. The
parameters of FPCA were set somewhat loosely for easy problems.
Specifically, we set $\bar\mu=10^{-4},xtol=10^{-4},\tau=2, I_m=10$
and all other parameters were set to the values given in Table
\ref{table:parameters-fpc}. Relative errors and times were averaged
over 5 runs. In this subsection, all test matrices were square,
i.e., $m=n.$

\begin{table}[ht]
\begin{center}\caption{Comparison of FPCA and SVT on easy problems}\label{table:Comparison-FPCA-SVT-easy}
\begin{tabular}{c c c c c | c c | c c}\hline

\multicolumn{5}{c|}{Problems} & \multicolumn{2}{c|}{FPCA} &
\multicolumn{2}{c}{SVT}
\\\hline

n & r & p &  SR & FR & rel.err. & time & rel.err. & time \\\hline

100 & 10 & 5666 & 0.57 & 0.34 & 4.27e-5 & 0.39 & 1.64e-3 &  30.40
\\\hline

200 & 10 & 15665 & 0.39 & 0.25 & 6.40e-5 & 1.38 & 1.90e-4 &  9.33
\\\hline

500 & 10 & 49471 & 0.20 & 0.20 & 2.48e-4 & 8.01 & 1.88e-4 &  23.77
\\\hline

1000 & 10 & 119406 & 0.12 & 0.17 & 5.04e-4 & 18.49 & 1.68e-4 & 41.81
\\\hline

1000 & 50 & 389852 & 0.39 & 0.25 & 3.13e-5 & 120.64 & 1.63e-4 &
228.79
\\\hline

1000 & 100 & 569900 & 0.57 & 0.33 & 2.26e-5 & 177.17 & 1.71e-4 &
635.15
\\\hline

5000 & 10 & 597973 & 0.02 & 0.17 & 1.58e-3 & 1037.12 & 1.73e-4 &
121.39\\\hline

5000 & 50 & 2486747 & 0.10 & 0.20 & 5.39e-4  &  1252.70 & 1.59e-4 &
1375.33
\\\hline

5000 & 100 & 3957533 & 0.16 & 0.25 & 2.90e-4 & 2347.41 & 1.74e-4 &
5369.76
\\\hline
\end{tabular}
\end{center}
\end{table}

From Table \ref{table:Comparison-FPCA-SVT-easy}, we can see that for
the easy problems except for one problem which is exceptionally
sparse as well as having low rank, FPCA was much faster and usually
provided more accurate solutions than SVT.

For hard problems, all parameters of FPCA were set to the values
given in Table \ref{table:parameters-fpc}, except that we set
$xtol=10^{-6}$ since this value is small enough to guarantee very
good recoverability. Also, for small problems ( i.e.,
$\max\{m,n\}<1000$ ), we set $I_m=500$; and for large problems (
i.e., $\max\{m,n\}\geq 1000$ ), we set $I_m=20.$ We use ``---'' to
indicate that the algorithm either diverges or does not terminate in
one hour. Relative errors and times were averaged over 5 runs.

\begin{table}[ht]
\begin{center}\caption{Comparison of FPCA and SVT on hard problems}\label{table:Comparison-FPCA-SVT-hard}
\begin{tabular}{c c c c | c c | c c}\hline

\multicolumn{4}{c|}{Problems} & \multicolumn{2}{c|}{FPCA} &
\multicolumn{2}{c}{SVT}
\\\hline

n & r & SR & FR & rel.err. & time & rel.err. & time \\\hline

40 & 9 & 0.5 & 0.80 & 1.21e-5 & 5.72 & 5.01e-1 &  3.05 \\\hline

100 & 14 & 0.3 & 0.87 & 1.32e-4 & 19.00 & 8.31e-1 &  316.90
\\\hline

1000 & 20 & 0.1 & 0.40 & 2.46e-5 & 116.15 & --- & --- \\\hline

1000 & 30 & 0.1 & 0.59 & 2.00e-3 & 128.30 & --- & --- \\\hline

1000 & 50 & 0.2 & 0.49 & 1.04e-5 & 183.67 & --- & --- \\\hline

\end{tabular}
\end{center}
\end{table}

From Table \ref{table:Comparison-FPCA-SVT-hard}, we can see that for
the hard problems, SVT either diverged or did not solve the problems
in less than one hour, or it yielded a very inaccurate solution. In
contrast, FPCA always provided a very good solution efficiently.

We can also see that FPCA was able to efficiently solve large
problems ($m=n=1000$) that could not be solved by SDPT3 due to the
large size of the matrices and the large number of constraints.

%
%
%
%
%
%
%
%
%
%
%

\subsection{Results for real data matrices}
In this section, we consider matrix completion problems based on two
real data sets: the Jester joke data set \citep{Goldberg-Roeder-Gupta-Perkins-2001} and the DNA data set \citep{Spellman-1998}. 
The Jester joke data set contains 4.1 million ratings for 100 jokes
from 73,421 users and is available on the website http://www.ieor.berkeley.edu/\~{}Egoldberg/jester-data/. 
Since the number of jokes is only 100, but the number of users is
quite large, we randomly selected $n_u$ users to get a modestly
sized matrix for testing purpose. As in
\citep{Srebro-Jaakkola-2003}, we randomly held out two ratings for
each user. Since some entries in the matrix are missing, we cannot
compute the relative error as we did for the randomly created
matrices. Instead, we computed the Normalized Mean Absolute Error
(NMAE) as in \citep{Goldberg-Roeder-Gupta-Perkins-2001} and
\citep{Srebro-Jaakkola-2003}. The Mean Absolute Error (MAE) is
defined as
\bea\label{eq:MAE-Jester}MAE=\df{1}{2N}\sum_{i=1}^N|\hat{r}_{i_1}^i-r_{i_1}^i|+|\hat{r}_{i_2}^i-r_{i_2}^i|,\eea
where $r_j^i$ and $\hat{r}_j^i$ are the withheld and predicted
ratings of movie $j$ by user $i$, respectively, for $j=i_1,i_2.$
NMAE is defined as
\bea\label{eq:NMAE-Jester}NMAE=\df{MAE}{r_{\max}-r_{\min}},\eea
where $r_{\min}$ and $r_{\max}$ are lower and upper bounds for the
ratings. Since all ratings are scaled to the range $[-10,+10]$, we
have $r_{\min}=-10$ and $r_{\max}=10.$

\begin{figure}
\centering \subfigure{
\includegraphics[scale=0.35]{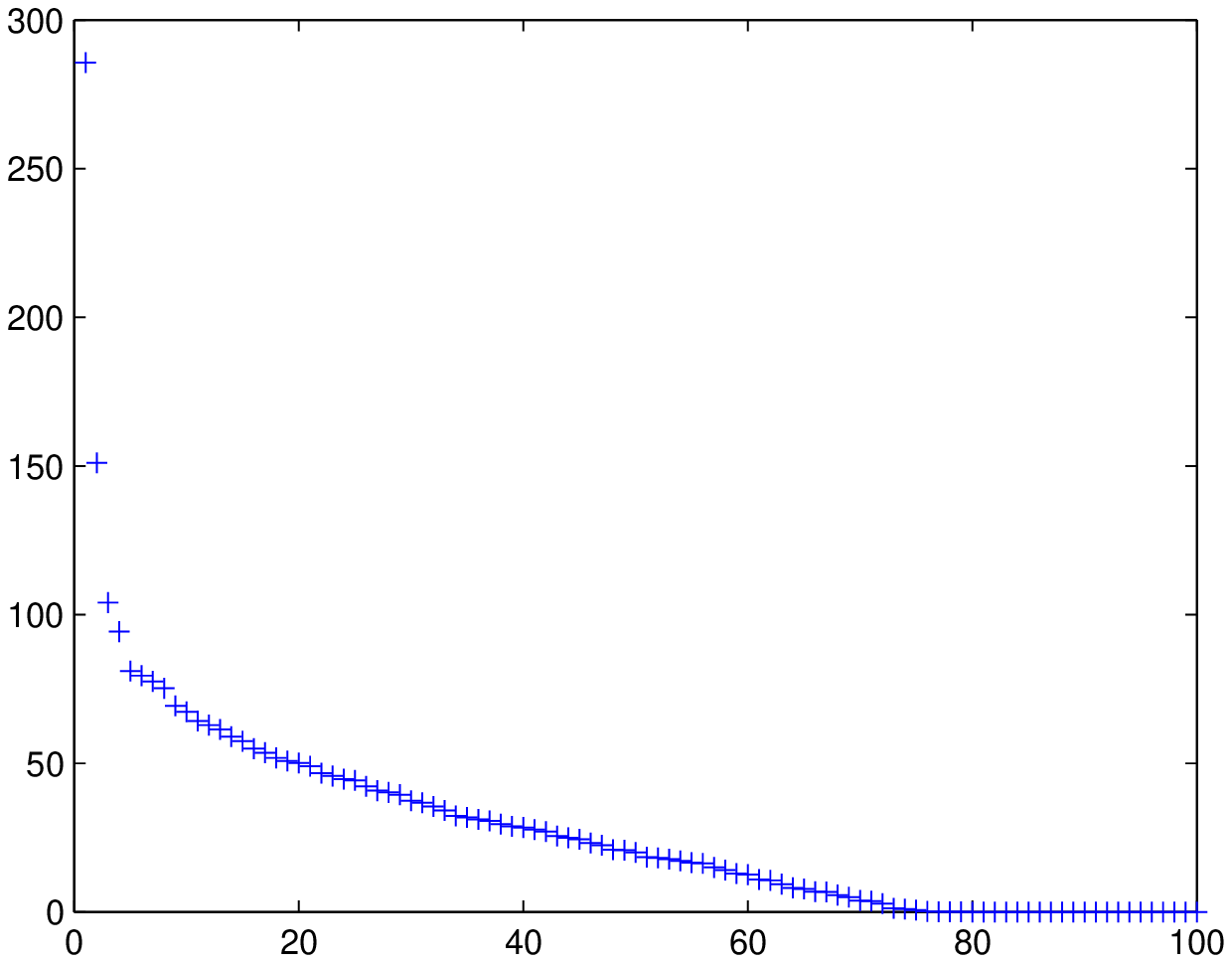}}\hspace{1mm}
\centering \subfigure{
\includegraphics[scale=0.35]{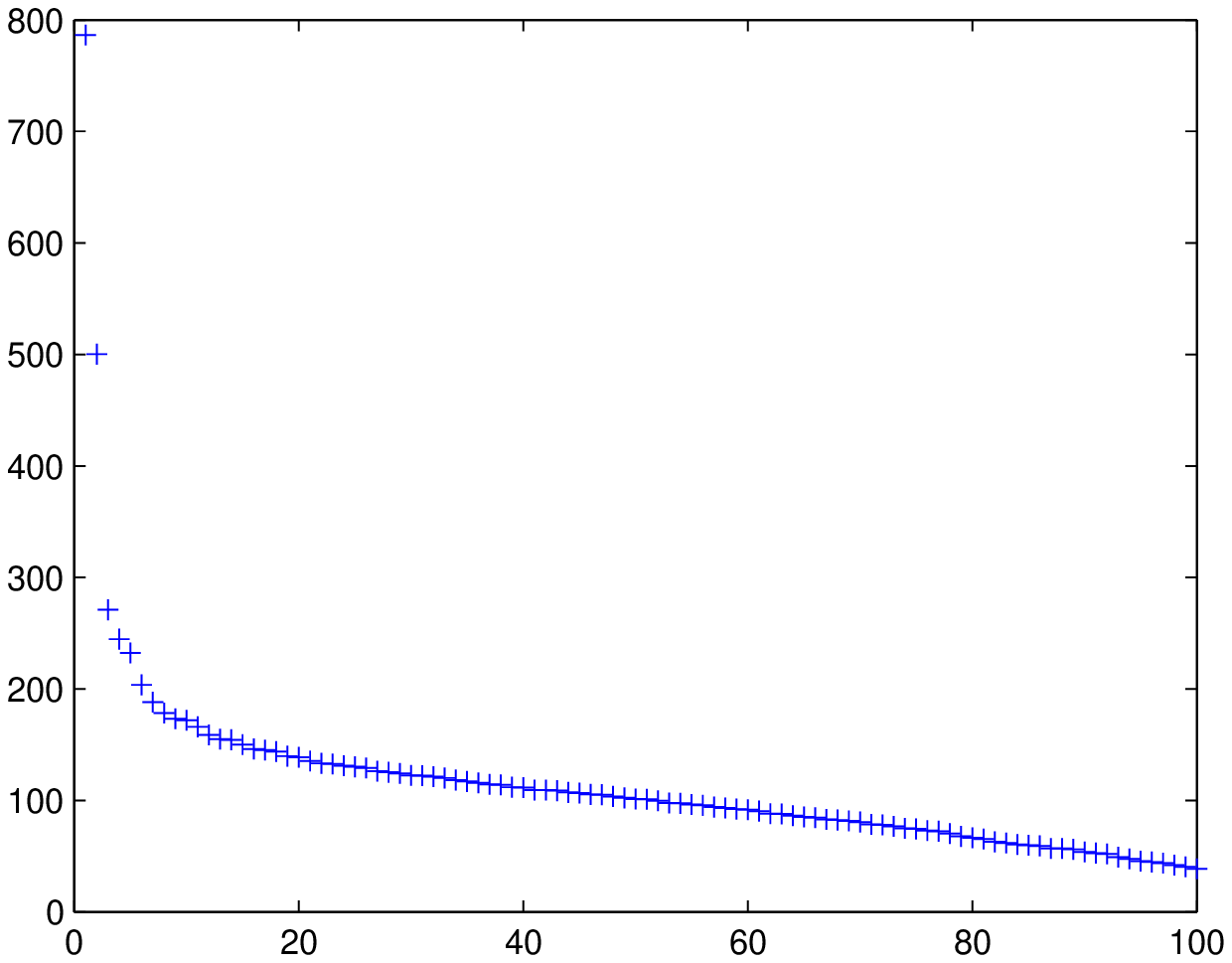}}\hspace{1mm}
\centering \subfigure{
\includegraphics[scale=0.35]{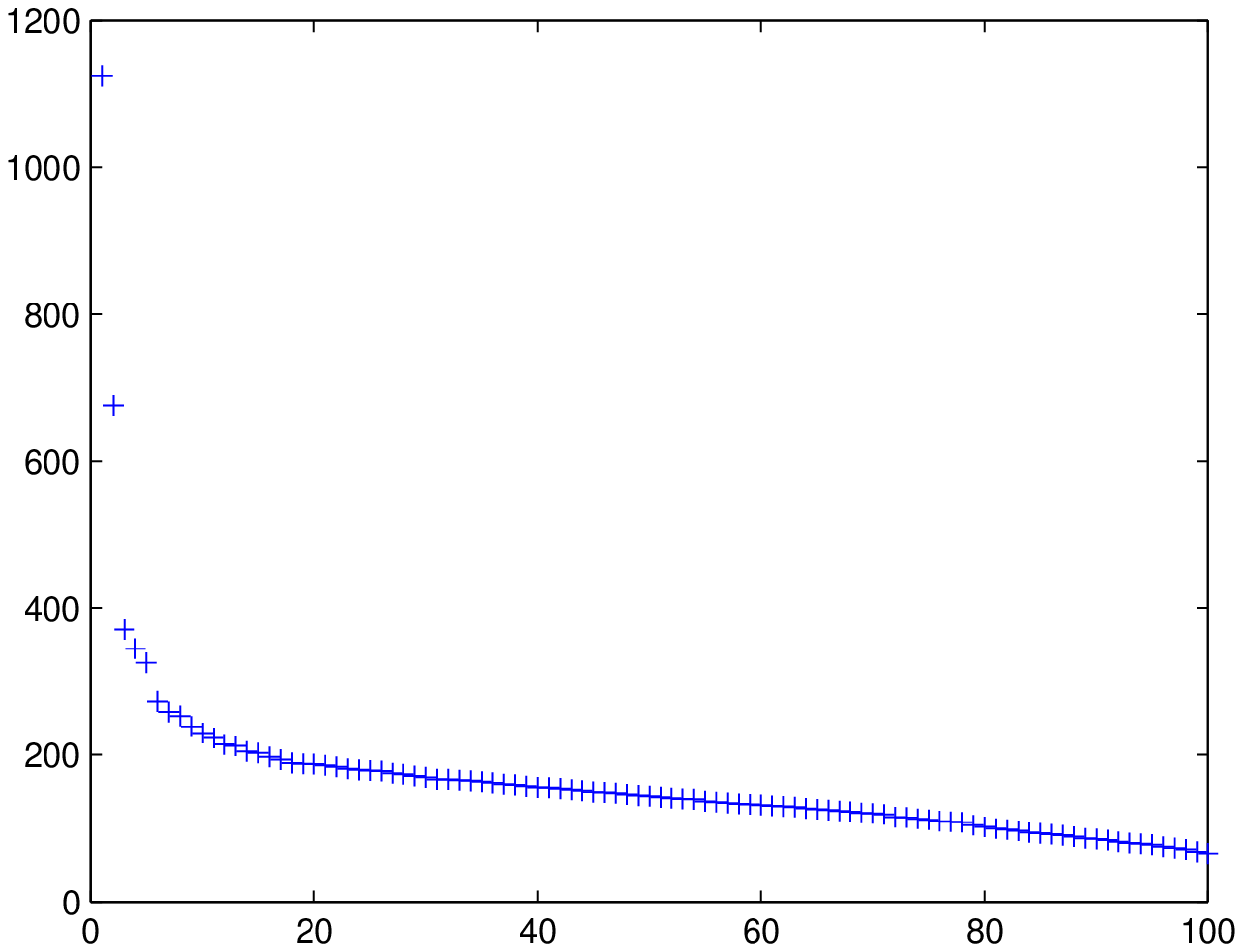}}
\caption{Distribution of the singular values of the recovered
matrices for the Jester data set using FPC1. Left:100 users, Middle:
1000 users, Right: 2000 users} \label{fig:jester_exact_svd}
\end{figure}

\begin{figure}
\centering \subfigure{
\includegraphics[scale=0.5]{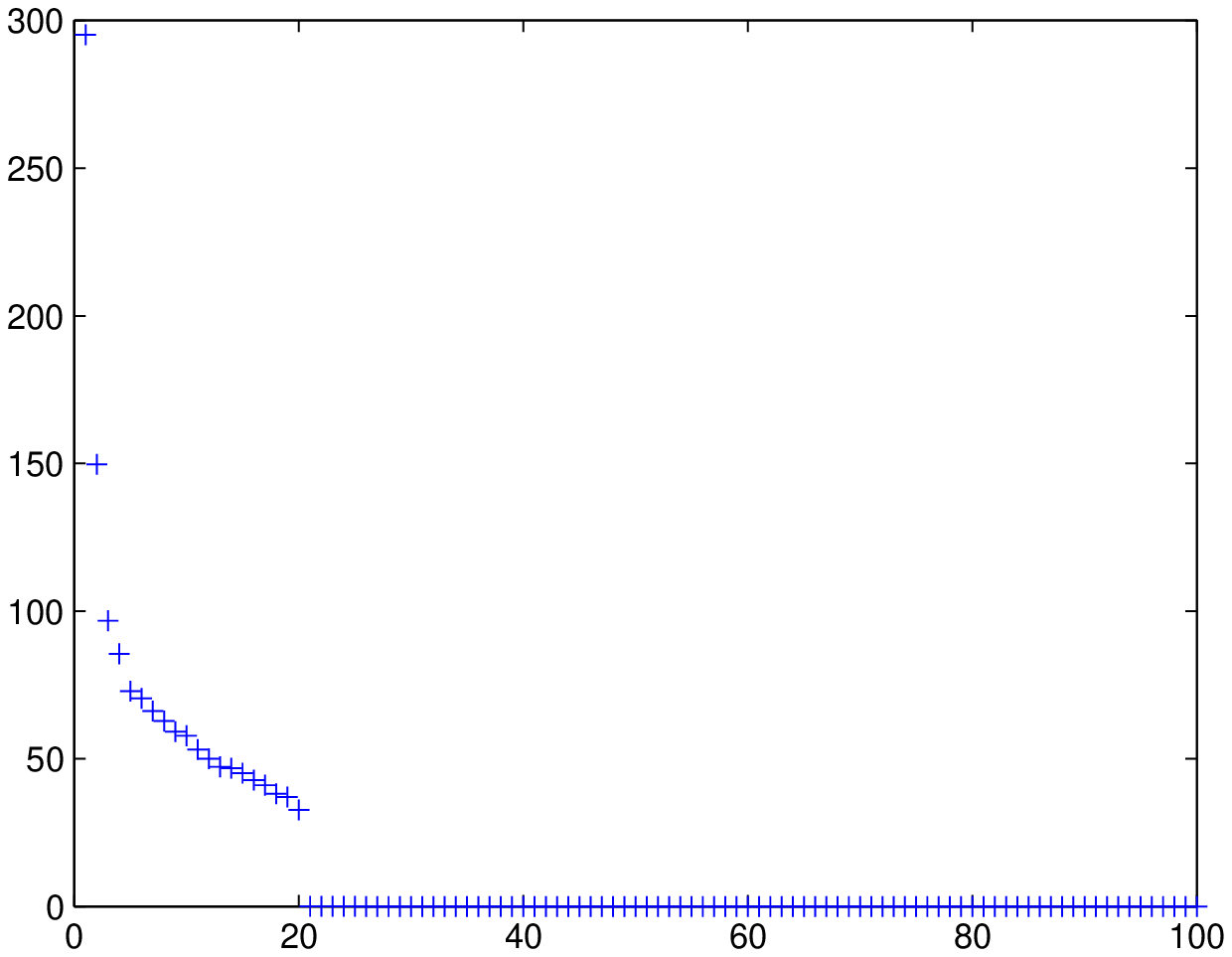}}
%
\subfigure{
\includegraphics[scale=0.5]{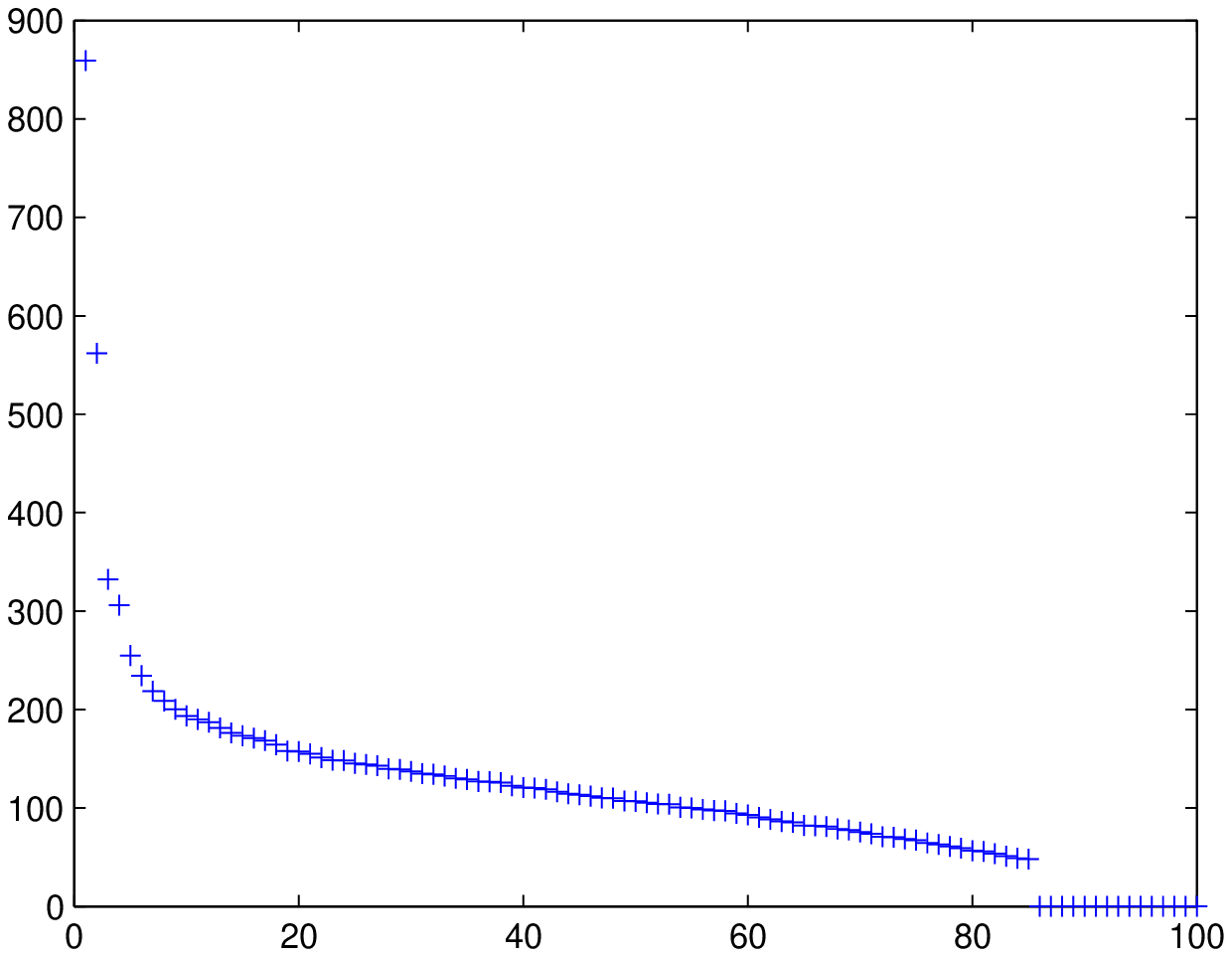}}
%
\subfigure{
\includegraphics[scale=0.5]{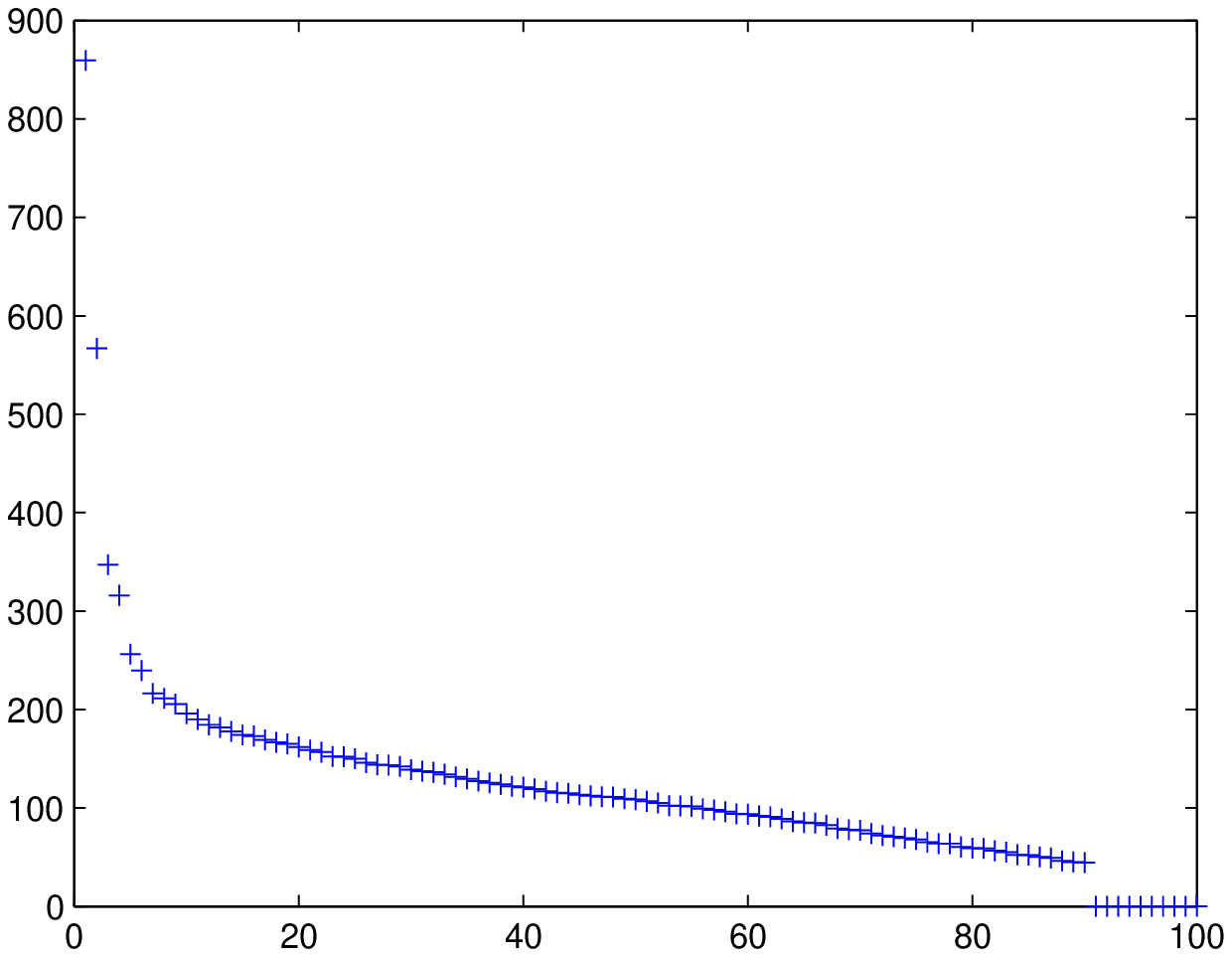}}
%
\subfigure{
\includegraphics[scale=0.5]{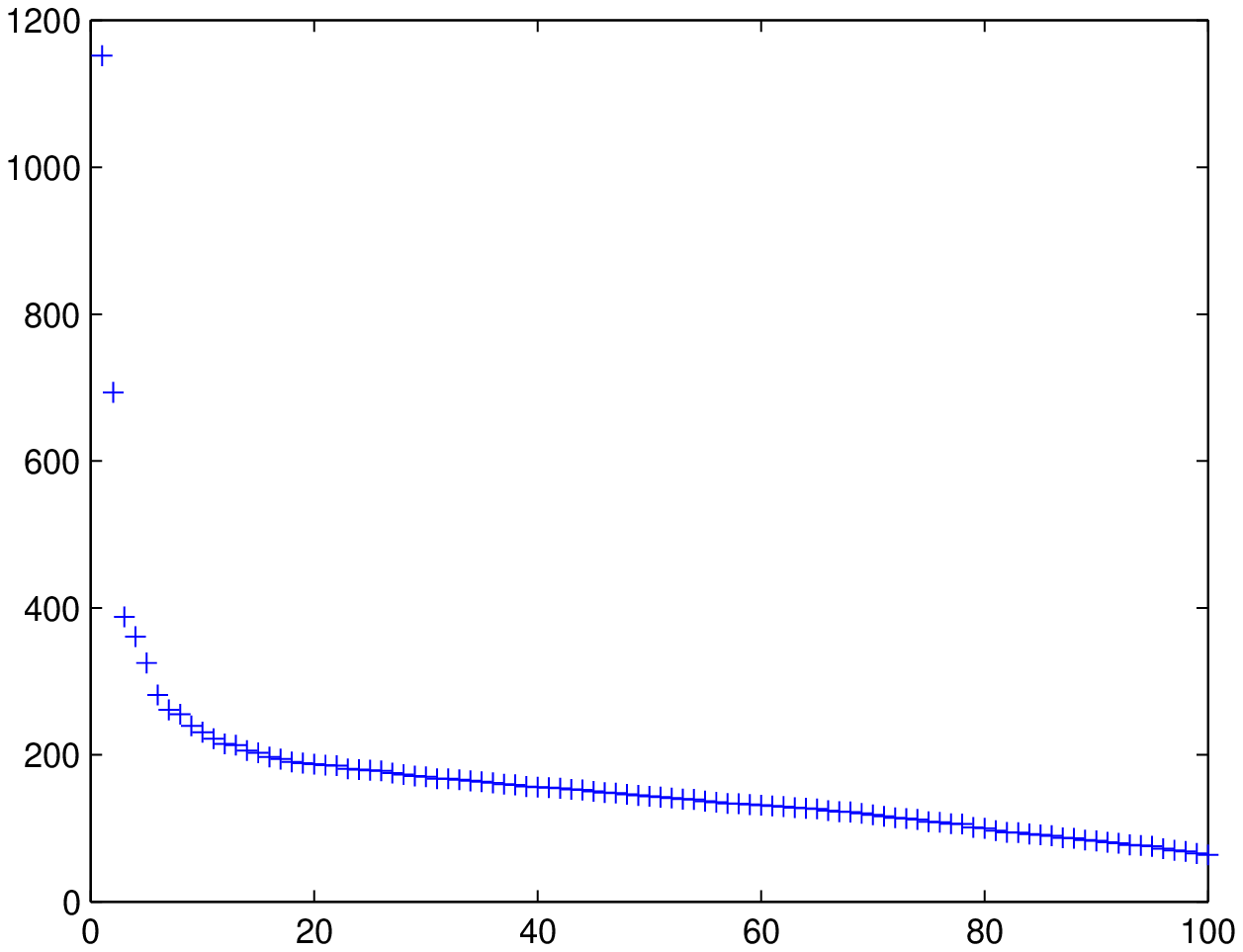}}
\caption{Distribution of the singular values of the recovered
matrices for the Jester data set using FPCA. Upper Left: 100 users,
$\epsilon_{k_s}=10^{-2},c_s=25$; Upper Right: 1000 users,
$\epsilon_{k_s}=10^{-2},c_s=100$; Bottom Left: 1000 users,
$\epsilon_{k_s}=10^{-4},c_s=100$; Bottom Right: 2000 users,
$\epsilon_{k_s}=10^{-4},c_s=200$} \label{fig:Jester_approx_svd}
\end{figure}

The numerical results for the Jester data set using FPC1 and FPCA
are presented in Tables \ref{table:FPC-Jester-exactSVD} and
\ref{table:FPC-Jester-fastsvd}, respectively. In these two tables,
$\sigma_{\max}$ and $\sigma_{\min}$ are the largest and smallest
positive singular values of the recovered matrices, and $rank$ is
the rank of the recovered matrices. The distributions of the
singular values of the recovered matrices are shown in Figures
\ref{fig:jester_exact_svd} and \ref{fig:Jester_approx_svd}. From
Tables \ref{table:FPC-Jester-exactSVD} and
\ref{table:FPC-Jester-fastsvd} we can see that by using FPC1 and
FPCA to recover these matrices, we can get relatively low NMAEs,
which are comparable to the results shown in
\citep{Srebro-Jaakkola-2003} and
\citep{Goldberg-Roeder-Gupta-Perkins-2001}.

\begin{table}[ht]
\begin{center}\caption{Numerical results for FPC1 for the Jester joke data set}\label{table:FPC-Jester-exactSVD}
\begin{tabular}{c c c c c c c c}\hline
num.user & num.samp & samp.ratio & rank & $\sigma_{\max}$ &
$\sigma_{\min}$ & NMAE & Time
\\\hline
100   & 7172 & 0.7172 & 79 & 285.65 & 3.49e-4 & 0.1727 & 34.30
\\

1000  & 71152 & 0.7115 & 100 & 786.37 & 38.43 & 0.1667 & 304.81\\

2000 & 140691 & 0.7035 & 100 & 1.1242e+3 & 65.06 & 0.1582 & 661.66
\\ \hline
\end{tabular}
\end{center}
\end{table}

\begin{table}[ht]
\begin{center}\caption{Numerical results for FPCA for the Jester joke data set ($c_s$ is the number of rows we picked for the approximate SVD)}\label{table:FPC-Jester-fastsvd}
\begin{tabular}{c c c c c c c c c c}\hline
num.user & num.samp & samp.ratio & $\epsilon_{k_s}$& $c_s$ & rank &
$\sigma_{\max}$ & $\sigma_{\min}$ & NMAE & Time
\\\hline

100 & 7172 & 0.7172 & $10^{-2}$ & 25 & 20 & 295.14 & 32.68 & 0.1627
&
26.73 \\

1000 & 71152 & 0.7115 & $10^{-2}$ & 100 & 85 & 859.27 & 48.04 & 0.2008 & 808.52 \\

1000 & 71152 & 0.7115 & $10^{-4}$ & 100 & 90 & 859.46 & 44.62 &
0.2101 &  778.56 \\


2000 & 140691 & 0.7035 & $10^{-4}$ & 200 & 100 & 1.1518e+3 & 63.52 &
0.1564 & 1.1345e+3 \\

\hline
\end{tabular}
\end{center}
\end{table}

\begin{figure}
\centering \subfigure{
\includegraphics[scale=0.5]{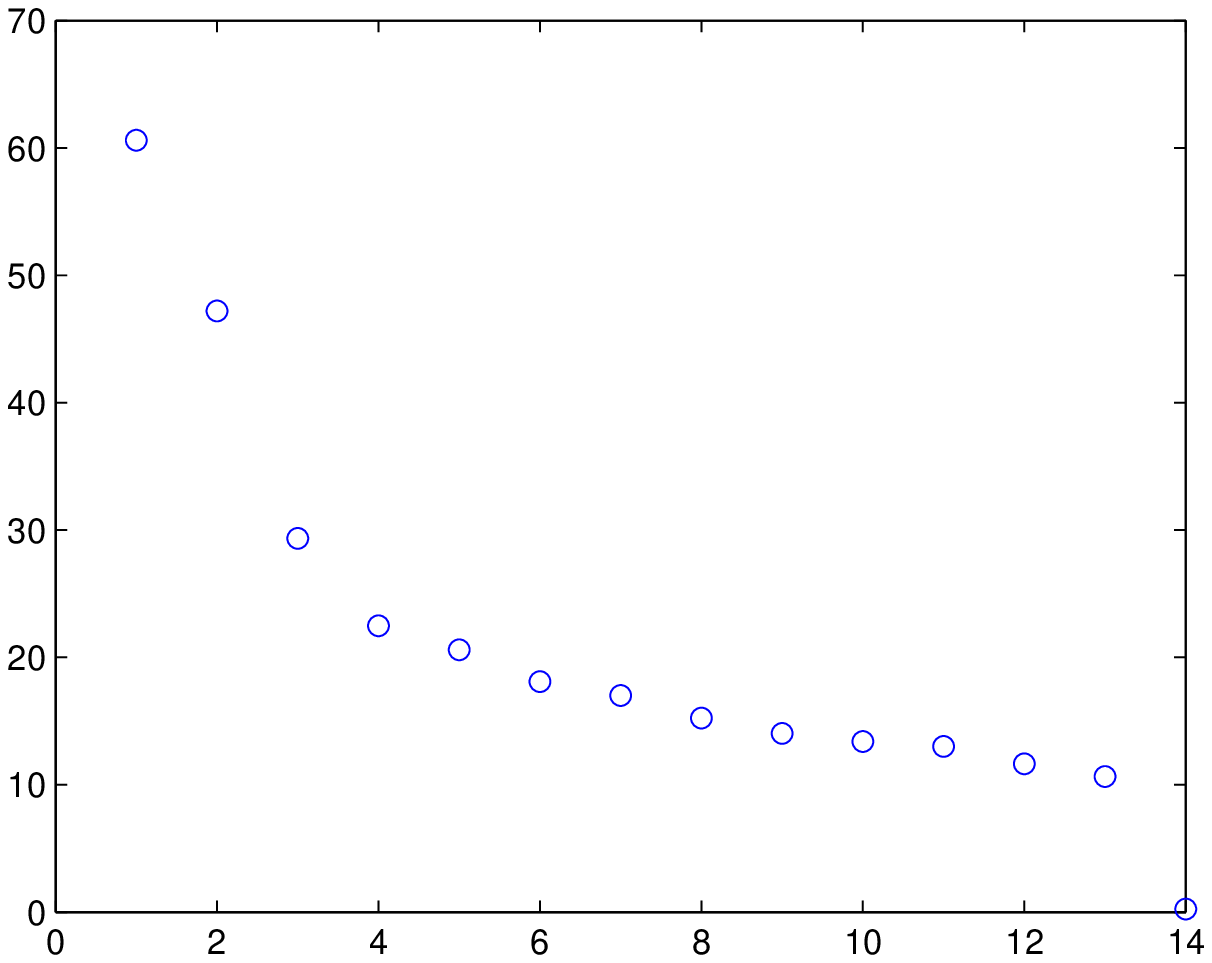}}
\subfigure{
\includegraphics[scale=0.5]{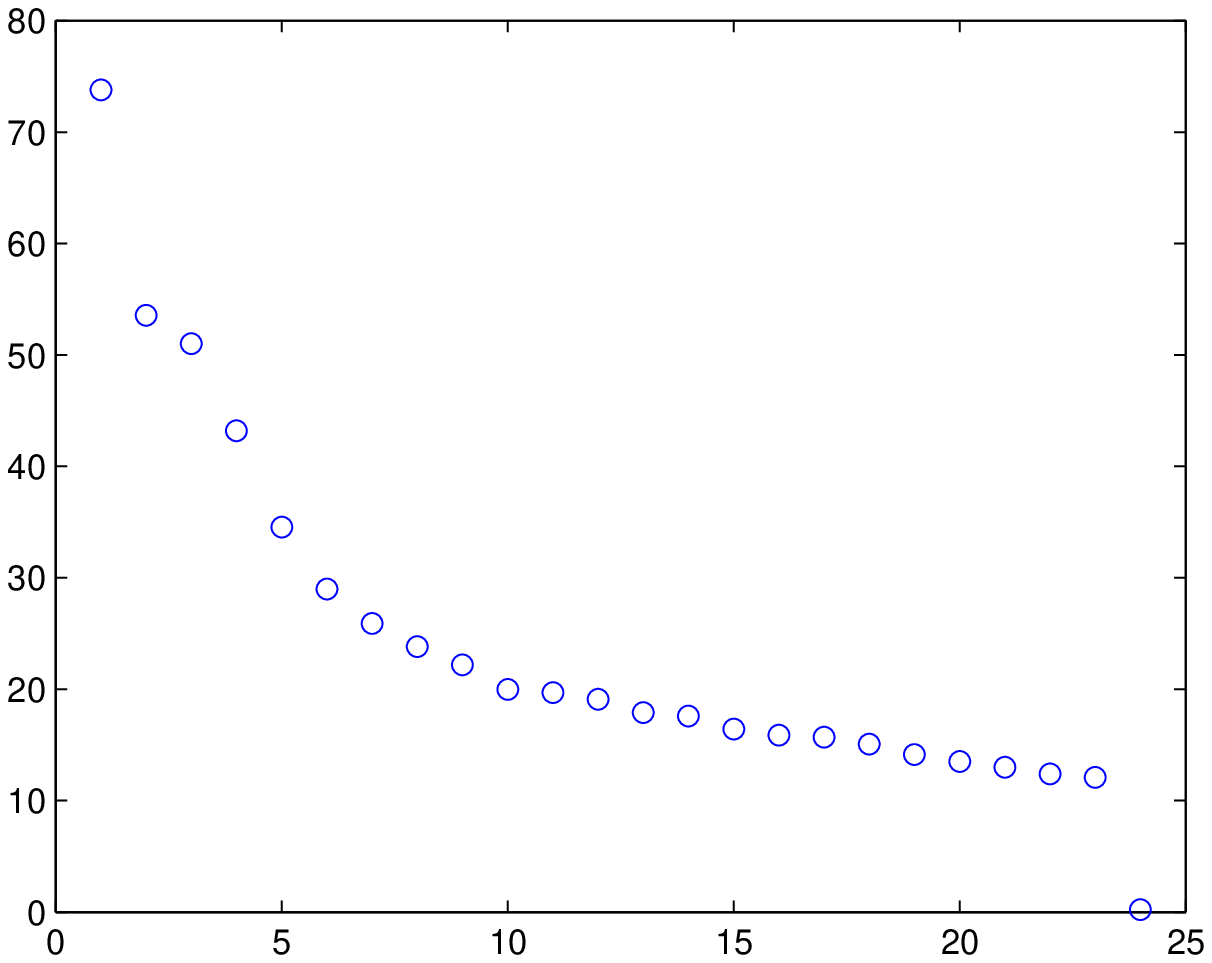}}
\caption{Distribution of the singular values of the matrices in the
original DNA microarray data sets. Left: Elutriation matrix; Right:
Cdc15 matrix. } \label{fig:DNA-actual-singular-values}
\end{figure}

We also used two data sets of DNA microarrays from
\citep{Spellman-1998}. These data sets are available on the website
http://cellcycle-www.stanford.edu/. The first microarray data set is
a matrix that represents the expression of 6178 genes in 14
experiments based on the Elutriation data set in
\citep{Spellman-1998}. The second microarray data set is based on
the Cdc15 data set in \citep{Spellman-1998}, and represents the
expression of 6178 genes in 24 experiments. However, some entries in
these two matrices are missing. For evaluating our algorithms, we
created complete matrices by deleted all rows containing missing
values. This is similar to how the DNA microarray data set was
preprocessed in \citep{Troyanskaya-2001}. The resulting complete
matrix for the Elutriation data set was $5766 \times 14$. The
complete matrix for the Cdc15 data set was $4381\times 24$. We must
point out that these DNA microarray matrices are neither low-rank
nor even approximately low-rank although such claims have been made
in some papers. The distributions of the singular values of these
two matrices are shown in Figure
\ref{fig:DNA-actual-singular-values}. From this figure we can see
that in each microarray matrix, only one singular value is close to
zero, while the others are far away from zero. Thus there is no way
to claim that the rank of the Elutriation matrix is less than 13, or
the rank of the Cdc15 matrix is less than 23. Since these matrices
are not low-rank, we cannot expect our algorithms to recover these
matrices by sampling only a small portion of their entries. Thus we
needed to further modify the data sets to yield low-rank matrices.
Specifically, we used the best rank-2 approximation to the
Elutriation matrix as the new complete Elutriation matrix and the
best rank-5 approximation to the Cdc15
matrix as the new complete Cdc15 matrix. 
The numerical results for FPCA for recovering these two matrices are
presented in Table \ref{table:num-res-fpc-DNA}. In the FPCA
algorithm, we set $\epsilon_{k_s}=10^{-2}$ and $xtol=10^{-6}$. For
the Elutriation matrix, we set $c_s=115$ and for the Cdc15 matrix,
we set $c_s=88$. The observed entries were randomly sampled. From
Table \ref{table:num-res-fpc-DNA} we can see that by taking 60\% of
the entries of the matrices, our FPCA algorithm can recover these
matrices very well, yielding relative errors as low as $10^{-5}$ and
$10^{-6}$, which is promising for practical use.

\begin{table}[ht]
\begin{center}\caption{Numerical results of FPCA for DNA microarray data sets}\label{table:num-res-fpc-DNA}
\begin{tabular}{c c c c c c c c c}\hline
Matrix & m & n & p & rank & SR & FR  & rel.err & Time
\\\hline
Elutriation & 5766 & 14 & 48434 & 2 & 0.6 & 0.2386 & 1.83e-5 &
218.01
\\

Cdc15 & 4381 & 24 & 63086 & 5 & 0.6 & 0.3487 & 7.95e-6 & 189.32 \\
\hline
\end{tabular}
\end{center}
\end{table}

\begin{figure}
\centering
\begin{tabular}{cc}
\begin{minipage}{3.0in}
\centering\includegraphics[scale=0.5]{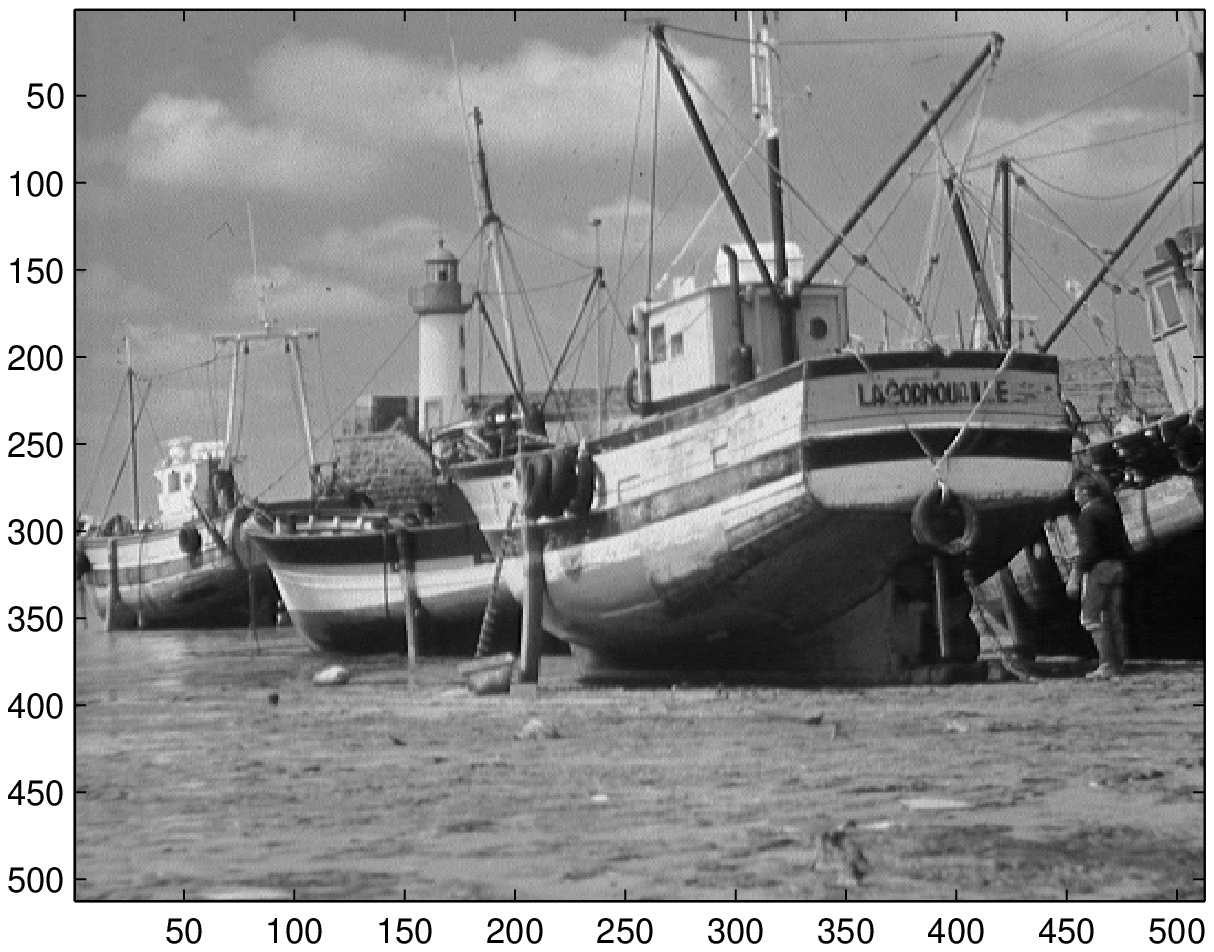}
\end{minipage}
&
\begin{minipage}{3.0in}
\centering\includegraphics[scale=0.5]{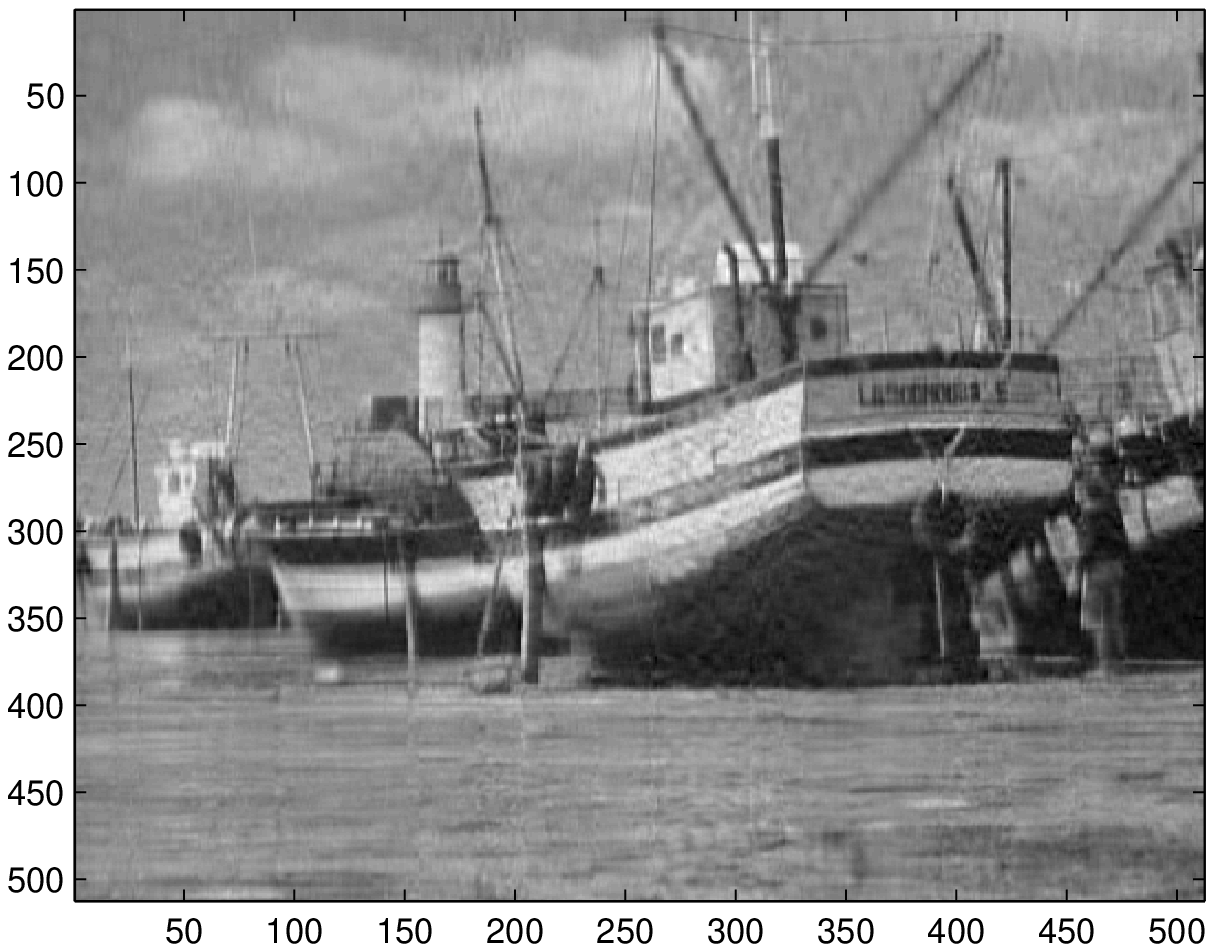}
\end{minipage}
\\
(a) & (b) \\
\begin{minipage}{3.0in}
\centering\includegraphics[scale=0.5]{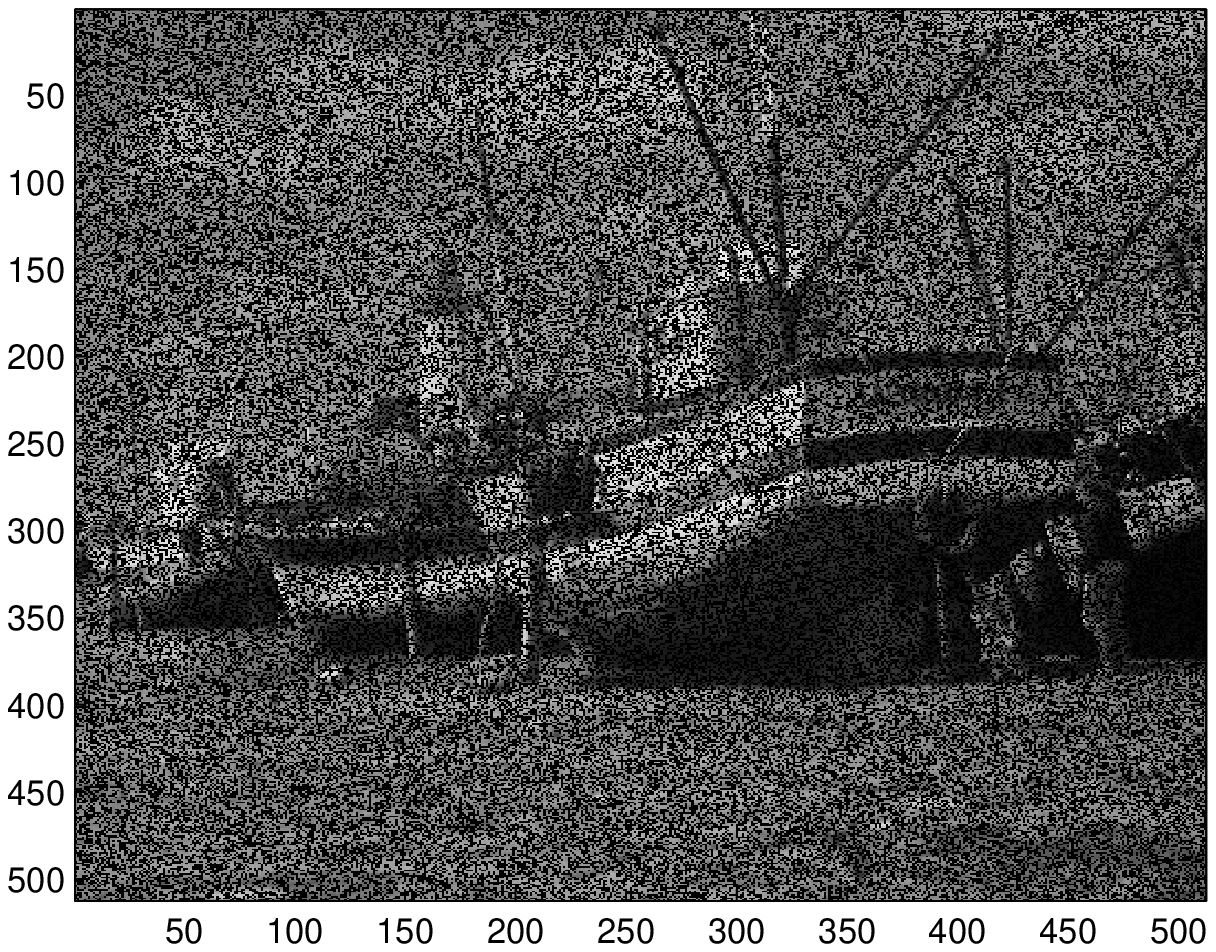}
\end{minipage}
&
\begin{minipage}{3.0in}
\centering\includegraphics[scale=0.5]{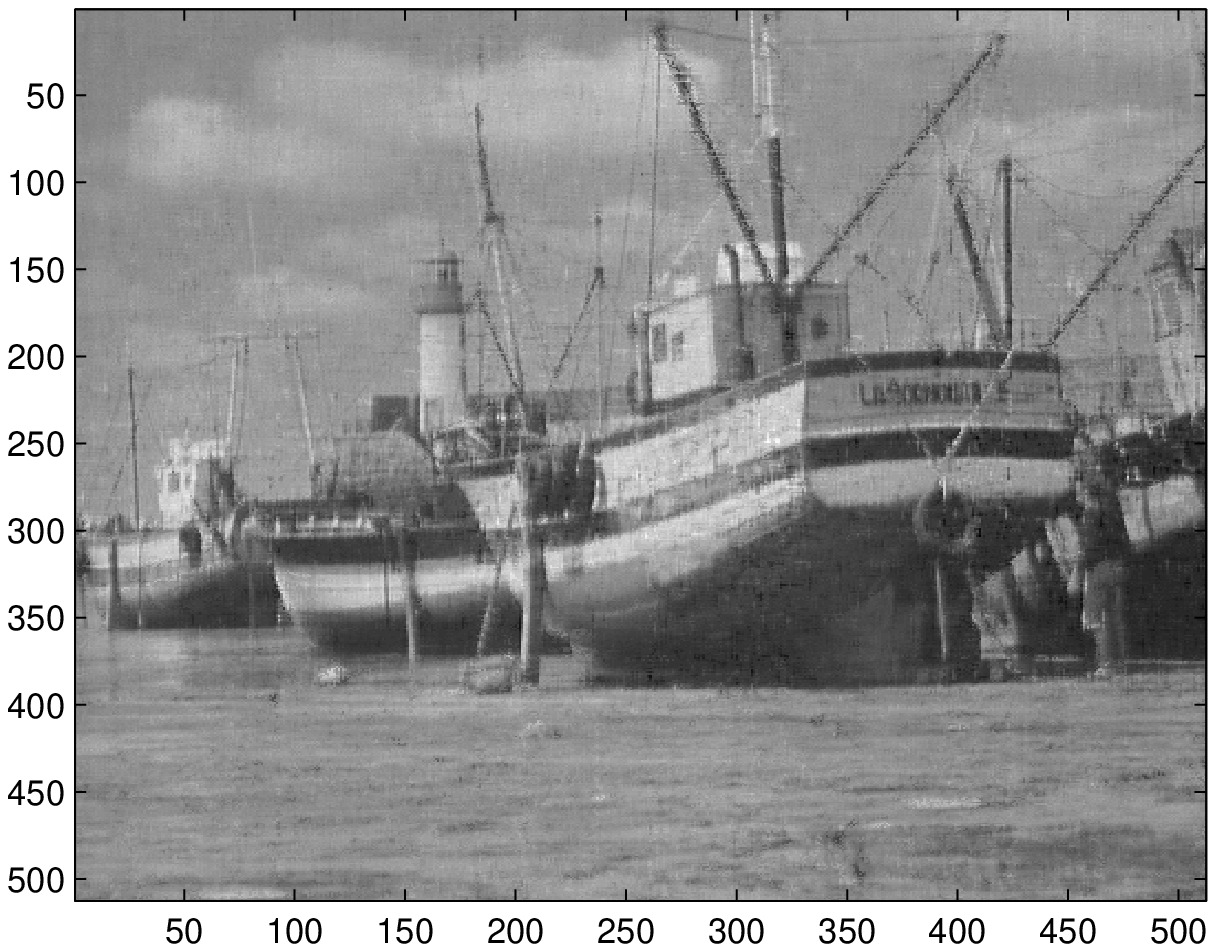}
\end{minipage}
\\
(c) & (d) \\
\begin{minipage}{3.0in}
\centering\includegraphics[scale=0.5]{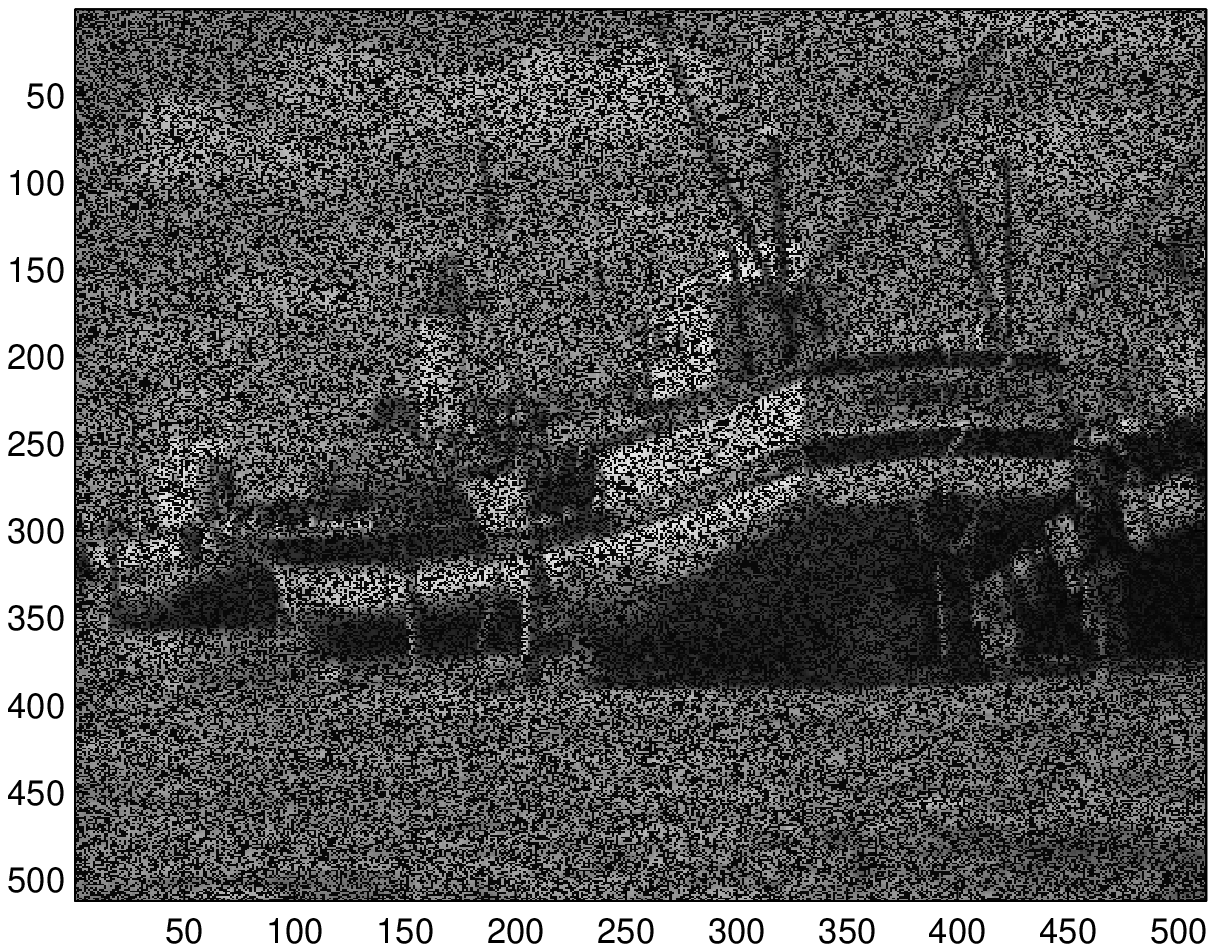}
\end{minipage}
&
\begin{minipage}{3.0in}
\centering\includegraphics[scale=0.5]{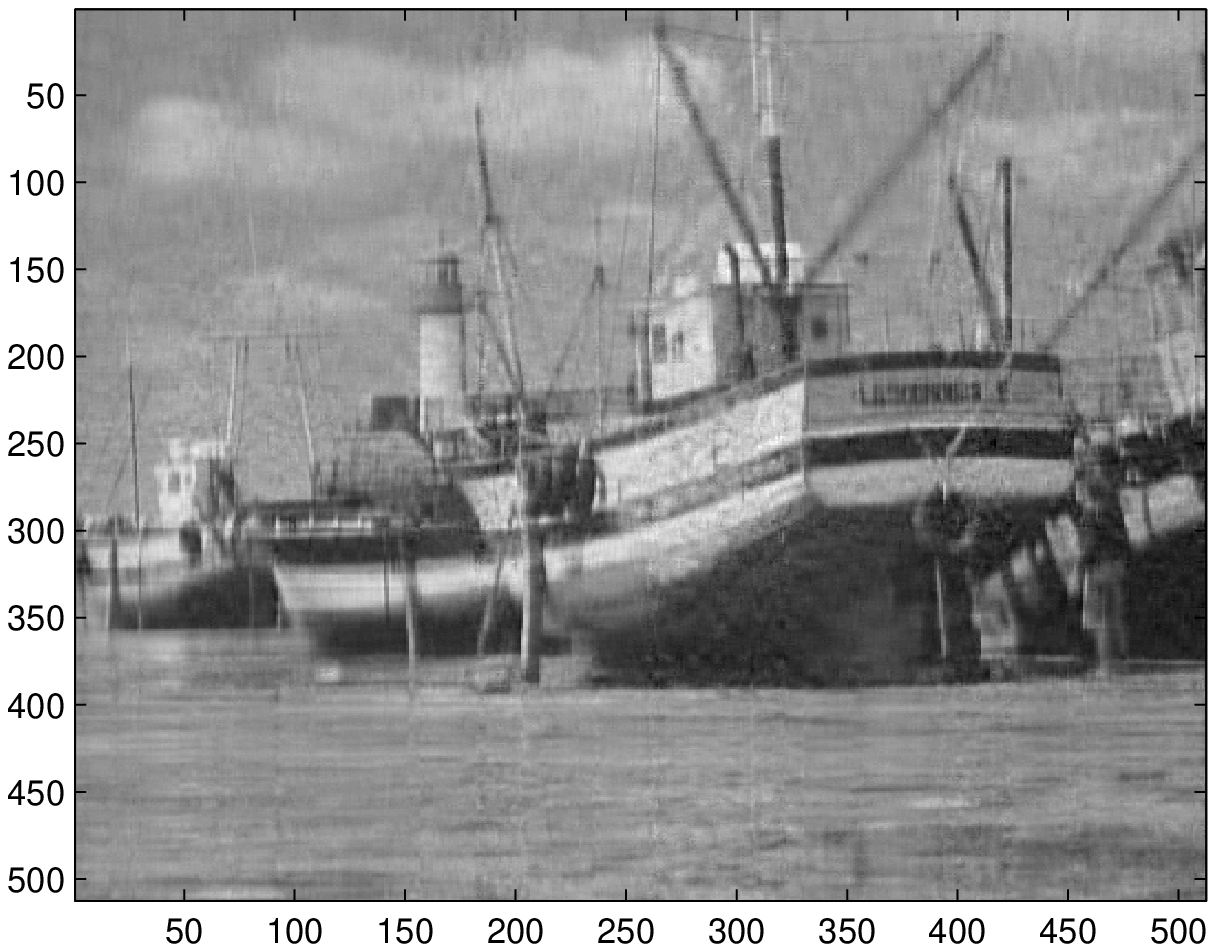}
\end{minipage}
\\
(e) & (f) \\
\begin{minipage}{3.0in}
\centering\includegraphics[scale=0.5]{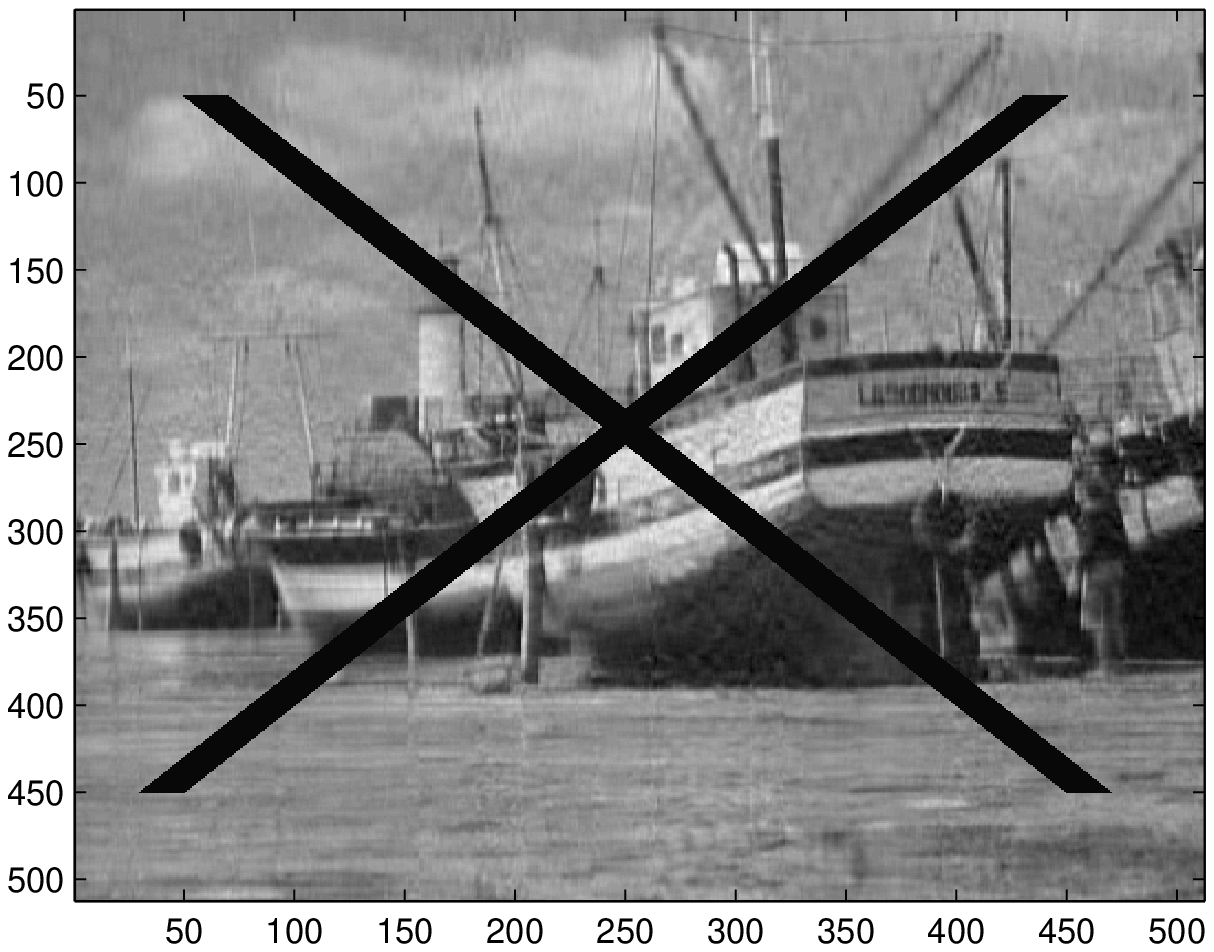}
\end{minipage}
&
\begin{minipage}{3.0in}
\centering\includegraphics[scale=0.5]{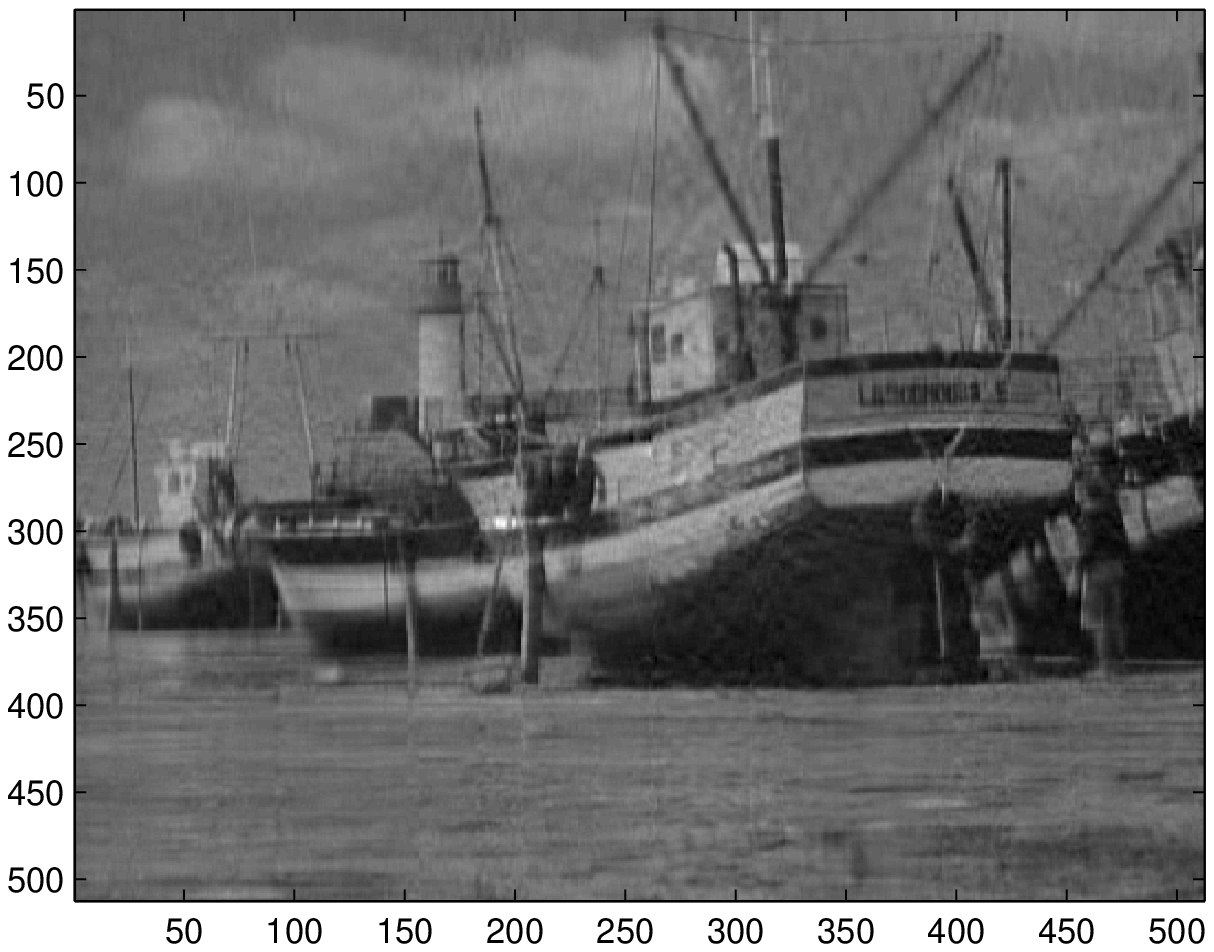}
\end{minipage}
\\
(g) & (h) \\
\end{tabular}
\caption{(a): Original $512\times 512$ image with full rank; (b):
Original image truncated to be of rank 40; (c): 50\% randomly masked
original image; (d): Recovered image from 50\% randomly masked
original image ($rel.err =8.41e-2$); (e): 50\% randomly masked rank
40 image; (f): Recovered image from 50\% randomly masked rank 40
image ($rel.err =3.61e-2$); (g): Deterministically masked rank 40
image (SR = 0.96); (h): Recovered image from deterministically
masked rank 40 image ($rel.err = 1.70e-2$).}\label{fig:boat}
\end{figure}

To graphically illustrate the effectiveness of FPCA, we applied it
to image inpainting
\citep{Bertalmio-Sapiro-Caselles-Ballester-2000}. Grayscale images
and color images can be expressed as matrices and tensors,
respectively. In grayscale image inpainting, the grayscale value of
some of the pixels of the image are missing, and we want to fill in
these missing values. If the image is of low-rank, or of numerical
low-rank, we can solve the image inpainting problem as a matrix
completion problem \eqref{prob:rank-matrix-completion-intro}. In our
test we applied SVD to the $512\times 512$ image in Figure
\ref{fig:boat}(a), and truncated this decomposition to get the
rank-40 image which is shown in Figure \ref{fig:boat}(b). Figure
\ref{fig:boat}(c) is a masked version of the image in Figure
\ref{fig:boat}(a), where one half of the pixels in Figure
\ref{fig:boat}(a) were masked uniformly at random. Figure
\ref{fig:boat}(d) is the image obtained from Figure
\ref{fig:boat}(c) by applying FPCA. Figure \ref{fig:boat}(d) is a
low-rank approximation to Figure \ref{fig:boat}(a) with a relative
error of $8.41e-2.$ Figure \ref{fig:boat}(e) is a masked version of
the image in Figure \ref{fig:boat}(b), where one half of the pixels
in Figure \ref{fig:boat}(b) were masked uniformly at random. Figure
\ref{fig:boat}(f) is the image obtained from Figure
\ref{fig:boat}(e) by applying FPCA. Figure \ref{fig:boat}(f) is an
approximation to Figure \ref{fig:boat}(b) with a relative error of
$3.61e-2.$ Figure \ref{fig:boat}(g) is another masked image obtained
from Figure \ref{fig:boat}(b), where 4 percent of the pixels were
masked in a non-random fashion. Figure \ref{fig:boat}(h) is the
image obtained from Figure \ref{fig:boat}(g) by applying FPCA.
Figure \ref{fig:boat}(g) is an approximation to Figure
\ref{fig:boat}(b) with a relative error of $1.70e-2$.


\section{Conclusions and discussions} In this paper, we derived a
fixed point continuation algorithm and a Bregman iterative algorithm
for solving the linearly constrained nuclear norm minimization
problem, which is a convex relaxation of the NP-hard linearly
constrained matrix rank minimization problem. The convergence of the
fixed point iterative scheme was established. By adopting an
approximate SVD technique, we obtained a very powerful algorithm
(FPCA) for the matrix rank minimization problem. On matrix
completion problems, FPCA greatly outperforms SDP solvers such as
SDPT3 in both speed and recoverability of low-rank matrices. Further
study is needed to prove the convergence of algorithm FPCA.

\acknowledgement We would like to thank two anonymous referees for
their helpful comments. The first author thanks Junzhou Huang from
Rutgers University for fruitful discussions on image inpainting.

\bibliographystyle{spmpsci}
\bibliography{matrix-rank-minimization}

\begin{thebibliography}{10}
\providecommand{\url}[1]{{#1}}
\providecommand{\urlprefix}{URL }
\expandafter\ifx\csname urlstyle\endcsname\relax
  \providecommand{\doi}[1]{DOI~\discretionary{}{}{}#1}\else
  \providecommand{\doi}{DOI~\discretionary{}{}{}\begingroup
  \urlstyle{rm}\Url}\fi

\bibitem{Bach-2008}
Bach, F.R.: Consistency of trace norm minimization.
\newblock Journal of Machine Learning Research \textbf{9}(Jun), 1019--1048
  (2008)

\bibitem{vandenBerg-Friedlander-2008}
van~den Berg, E., Friedlander, M.P.: Probing the {P}areto frontier for basis
  pursuit solutions.
\newblock Preprint available at Optimization Online: 2008.01.1889  (2008)

\bibitem{Bertalmio-Sapiro-Caselles-Ballester-2000}
Bertalm\'io, M., Sapiro, G., Caselles, V., Ballester, C.: Image inpainting.
\newblock Proceedings of SIGGRAPH 2000, New Orleans, USA  (2000)

\bibitem{Borwein-Lewis-2000-book}
Borwein, J.M., Lewis, A.S.: Convex Analysis and Nonlinear Optimization.
\newblock Springer-Verlag (2003)

\bibitem{Bregman-1967}
Bregman, L.: The relaxation method of finding the common points of convex sets
  and its application to the solution of problems in convex programming.
\newblock USSR Computational Mathematics and Mathematical Physics \textbf{7},
  200--217 (1967)

\bibitem{Burer-Monteiro-2003}
Burer, S., Monteiro, R.D.C.: A nonlinear programming algorithm for solving
  semidefinite programs via low-rank factorization.
\newblock Mathematical Programming (Series B) \textbf{95}, 329--357 (2003)

\bibitem{Burer-Monteiro-2005}
Burer, S., Monteiro, R.D.C.: Local mimima and convergence in low-rank
  semidefinite programming.
\newblock Mathematical Programming \textbf{103}(3), 427--444 (2005)

\bibitem{Cai-Candes-Shen-2008}
Cai, J., Cand\`es, E.J., Shen, Z.: A singular value thresholding algorithm for
  matrix completion.
\newblock submitted for publication  (2008)

\bibitem{Candes-Recht-2008}
Cand\`es, E.J., Recht, B.: Exact matrix completion via convex optimization.
\newblock Submitted  (2008)

\bibitem{Candes-Romberg-2005-l1-magic}
Cand\`es, E.J., Romberg, J.: $\ell_1$-{MAGIC}: Recovery of sparse signals via
  convex programming.
\newblock Tech. rep., Caltech (2005)

\bibitem{Candes-Romberg-Tao-2006}
Cand\`es, E.J., Romberg, J., Tao, T.: Robust uncertainty principles: {E}xact
  signal reconstruction from highly incomplete frequency information.
\newblock IEEE Transactions on Information Theory \textbf{52}, 489--509 (2006)

\bibitem{Candes-Tao-2009}
Cand\`es, E.J., Tao, T.: The power of convex relaxation: near-optimal matrix
  completion.
\newblock preprint  (2009)

\bibitem{Dai-Milenkovie-2008}
Dai, W., Milenkovic, O.: Subspace pursuit for compressive sensing: closing the
  gap between performance and complexity.
\newblock Preprint available at arXiv: 0803.0811  (2008)

\bibitem{Donoho-2006}
Donoho, D.: Compressed sensing.
\newblock IEEE Transactions on Information Theory \textbf{52}, 1289--1306
  (2006)

\bibitem{Donoho-Tsaig-Drori-Starck-2006}
Donoho, D., Tsaig, Y., Drori, I., Starck, J.C.: Sparse solution of
  underdetermined linear equations by stagewise orthogonal matching pursuit.
\newblock Submitted to IEEE Trransactions on Information Theory  (2006)

\bibitem{Donoho-Tsaig-2006}
Donoho, D.L., Tsaig, Y.: Fast solution of $\ell_1$-norm minimization problems
  when the solution may be sparse.
\newblock Tech. rep., Department of Statistics, Stanford University (2006)

\bibitem{Drineas-Kannan-Mahoney-2006}
Drineas, P., Kannan, R., Mahoney, M.W.: Fast {M}onte {C}arlo algorithms for
  matrices ii: Computing low-rank approximations to a matrix.
\newblock SIAM J. Computing \textbf{36}, 132--157 (2006)

\bibitem{Fazel-thesis-2002}
Fazel, M.: Matrix rank minimization with applications.
\newblock Ph.D. thesis, Stanford University (2002)

\bibitem{Fazel-Hindi-Boyd-2001}
Fazel, M., Hindi, H., Boyd, S.: A rank minimization heuristic with application
  to minimum order system approximation.
\newblock In: Proceedings of the American Control Conference (2001)

\bibitem{Figueiredo-Nowak-Wright-2007}
Figueiredo, M.A.T., Nowak, R.D., Wright, S.J.: Gradient projection for sparse
  reconstruction: Application to compressed sensing and other inverse problems.
\newblock IEEE Journal on Selected Topics in Signal Processing \textbf{1}(4)
  (2007)

\bibitem{ElGhaoui-Gahinet-1993}
Ghaoui, L.E., Gahinet, P.: Rank minimization under {LMI} constraints: A
  framework for output feedback problems.
\newblock In: Proceedings of the European Control Conference (1993)

\bibitem{Goldberg-Roeder-Gupta-Perkins-2001}
Goldberg, K., Roeder, T., Gupta, D., Perkins, C.: Eigentaste: A constant time
  collaborative filtering algorithm.
\newblock Information Retrieval \textbf{4}(2), 133--151 (2001)

\bibitem{Goldfarb-Ma-2009}
Goldfarb, D., Ma, S.: Convergence of fixed point continuation algorithms for
  matrix rank minimization.
\newblock Tech. rep., Department of IEOR, Columbia University (2009)

\bibitem{Hale-Yin-Zhang-2007}
Hale, E.T., Yin, W., Zhang, Y.: A fixed-point continuation method for
  $\ell_1$-regularized minimization with applications to compressed sensing.
\newblock Tech. rep., CAAM TR07-07 (2007)

\bibitem{Hiriart-Urruty-Lemarechal-1993}
Hiriart-Urruty, J.B., Lemar\'echal, C.: Convex Analysis and Minimization
  Algorithms II: Advanced Theory and Bundle Methods.
\newblock Springer-Verlag, New York (1993)

\bibitem{Horn-Johnson-book-1985}
Horn, R.A., Johnson, C.R.: Matrix Analysis.
\newblock Cambridge University Press (1985)

\bibitem{Kim-Koh-Lustig-Boyd-Gorinevsky-2007}
Kim, S.J., Koh, K., Lustig, M., Boyd, S., Gorinevsky, D.: A method for
  large-scale $\ell_1$-regularized least-squares.
\newblock IEEE Journal on Selected Topics in Signal Processing \textbf{4}(1),
  606--617 (2007)

\bibitem{Linial-London-Rabinovich-1995}
Linial, N., London, E., Rabinovich, Y.: The geometry of graphs and some of its
  algorithmic applications.
\newblock Combinatorica \textbf{15}, 215--245 (1995)

\bibitem{Liu-Vandenberghe-2008}
Liu, Z., Vandenberghe, L.: Interior-point method for nuclear norm approximation
  with application to system identification.
\newblock Submitted to Mathematical Programming Series B  (2008)

\bibitem{Natarajan-1995}
Natarajan, B.K.: Sparse approximation solutions to linear systems.
\newblock SIAM J. Computing \textbf{24}(2), 227--234 (1995)

\bibitem{Osher-Burger-Goldfarb-Xu-Yin-2005}
Osher, S., Burger, M., Goldfarb, D., Xu, J., Yin, W.: An iterative
  regularization method for total varitaion-based image restoration.
\newblock SIAM MMS \textbf{4}(2), 460--489 (2005)

\bibitem{Recht-Fazel-Parrilo-2007}
Recht, B., Fazel, M., Parrilo, P.: Guaranteed minimum rank solutions of matrix
  equations via nuclear norm minimization.
\newblock Submitted to SIAM Review  (2007)

\bibitem{Rennie-Srebro-2005}
Rennie, J.D.M., Srebro, N.: Fast maximum margin matrix factorization for
  collaborative prediction.
\newblock In: Proceedings of the International Conference of Machine Learning
  (2005)

\bibitem{Rudin-Osher-Fatemi-1992}
Rudin, L., Osher, S., Fatemi, E.: Nonlinear total variation based noise removal
  algorithms.
\newblock Physica D \textbf{60}, 259--268 (1992)

\bibitem{Spellman-1998}
Spellman, P.T., Sherlock, G., Zhang, M.Q., Iyer, V.R., Anders, K., Eisen, M.B.,
  Brown, P.O., Botstein, D., Futcher, B.: Comprehensive identification of cell
  cycle-regulated genes of the yeast saccharomyces cerevisiae by microarray
  hybridization.
\newblock Molecular Biology of the Cell \textbf{9}, 3273--3297 (1998)

\bibitem{Srebro-thesis-2004}
Srebro, N.: Learning with matrix factorizations.
\newblock Ph.D. thesis, Massachusetts Institute of Technology (2004)

\bibitem{Srebro-Jaakkola-2003}
Srebro, N., Jaakkola, T.: Weighted low-rank approximations.
\newblock In: Proceedings of the Twentieth International Conference on Machine
  Learning (ICML-2003) (2003)

\bibitem{Sturm-1999}
Sturm, J.F.: Using {SeDuMi} 1.02, a {M}atlab toolbox for optimization over
  symmetric cones.
\newblock Optimization Methods and Software \textbf{11-12}, 625--653 (1999)

\bibitem{Tibshirani-1996}
Tibshirani, R.: Regression shrinkage and selection via the lasso.
\newblock Journal Royal Statistical Society B \textbf{58}, 267--288 (1996)

\bibitem{Tropp-2006}
Tropp, J.: Just relax: Convex programming methods for identifying sparse
  signals.
\newblock IEEE Transactions on Information Theory \textbf{51}, 1030--1051
  (2006)

\bibitem{Troyanskaya-2001}
Troyanskaya, O., Cantor, M., Sherlock, G., Brown, P., Hastie, T., Tibshirani,
  R., Botstein, D., Altman, R.B.: Missing value estimation methods for {DNA}
  microarrays.
\newblock Bioinformatics \textbf{17}(6), 520--525 (2001)

\bibitem{Tutuncu-Toh-Todd-2003}
T\"ut\"unc\"u, R.H., Toh, K.C., Todd, M.J.: Solving
  semidefinite-quadratic-linear programs using {SDPT}3.
\newblock Mathematical Programming Series B \textbf{95}, 189--217 (2003)

\bibitem{Wen-Yin-Goldfarb-Zhang-2009}
Wen, Z., Yin, W., Goldfarb, D., Zhang, Y.: A fast algorithm for sparse
  reconstruction based on shrinkage, subspace optimization and continuation.
\newblock Tech. rep., Department of IEOR, Columbia University (2009)

\bibitem{Yin-Osher-Goldfarb-Darbon-2008}
Yin, W., Osher, S., Goldfarb, D., Darbon, J.: Bregman iterative algorithms for
  $\ell_1$-minimization with applications to compressed sensing.
\newblock SIAM J. Imaging Sci \textbf{1}(1), 143--168 (2008)

\end{thebibliography}

\end{document}